\documentclass[11pt]{amsart}
\usepackage{epsf}

\newcommand{\Ext}{Ext}
\newcommand{\mph}{\vphantom{\Bigg|}}
\renewcommand{\phi}{\varphi}
\renewcommand{\epsilon}{\varepsilon}
\newcommand{\suml}{\sum\limits}
\newenvironment{enum}{\begin{description}}{\end{description}}
\newcommand{\wCP}{\widetilde{\CP}}
\newcommand{\wz}{\widetilde z}
\newcommand{\wf}{\widetilde f}
\newcommand{\hf}{\widehat f}
\newcommand{\CP}{{\C P}}
\newcommand{\GL}{G\!L}
\newcommand{\simto}{\overset{\sim}{\to}}
\newcommand{\PDif}{\Psi\Dif}
\newcommand{\nl}{\boldsymbol{\lambda}}
\newcommand{\DG}{{\mathrm{DG}}}
\newcommand{\DR}{{\mathrm{DR}}}
\newcommand{\super}{{\mathrm{super}}}
\newcommand{\sing}{{\mathrm{sing}}}
\newcommand{\nto}[1]{\overset{#1}{\to}}
\newcommand{\op}{{\overline\partial}}
\newcommand{\ndot}{\bullet}

\newcommand{\A}{{\mathfrak A}}
\newcommand{\g}{{\mathfrak g}}

\newcommand{\gl}{{\mathfrak{gl}}}

\newcommand{\hPsi}{\widehat\Psi}
\newcommand{\wa}{\widetilde a}

\def\O{\mathcal O}

\def\matho#1{\mathop{\mathrm{#1}}}

\def\suml{\sum\limits}
\DeclareMathSizes{11.1}{10}{8}{6}

\newcommand\nfrac[2]
{\dfrac{\raisebox{-1pt}{\fontsize{11.1}{10pt}\selectfont$#1$}}
       {\raisebox{ 1pt}{\fontsize{11.1}{10pt}\selectfont$#2$}}}

\newcommand{\Vect}{{\matho{Vect}}}
\newcommand{\Alt}{\matho{Alt}}
\newcommand{\Ad}{\matho{Ad}}
\newcommand{\Ker}{\matho{Ker}}

\newcommand{\Der}{\matho{Der}}
\newcommand{\Lie}{{\matho{Lie}}}

\newcommand{\Altl}{\Alt\limits}
\newcommand{\Tr}{{\matho{Tr}}}

\newcommand{\id}{\matho{id}}
\newcommand{\Id}{\matho{Id}}
\newcommand{\ad}{\matho{ad}}
\newcommand{\fin}{{\mathrm{fin}}}

\newcommand{\even}{\mathrm{even}}

\newcommand{\Circle}{\mathrm{Circle}}
\newcommand{\Dif}{\mathrm{Dif}}

\newcommand{\C}{\mathbb C}
\newcommand{\R}{\mathbb R}
\newcommand{\Z}{\mathbb Z}
\newcommand{\D}{\mathcal D}

\newtheorem*{theorem}{Theorem}
\newtheorem*{propos}{Proposition}
\newtheorem*{lemma}{Lemma}
\newtheorem*{corollary}{Corollary}
\newtheorem*{conjecture}{Conjecture}

\theoremstyle{remark}
\newtheorem*{remark}{Remark}
\newtheorem*{example}{Example}

\theoremstyle{definition}
\newtheorem*{defin}{Definition}

\author{Boris Shoikhet}
\title[Integration of the Lifting formulas]%
{Integration of the Lifting formulas and the cyclic homology of the
algebras of differential operators}
\date{September 1, 1998}
\address{IUM, 11 Bol'shoj Vlas'evskij per.,
Moscow 121002, Russia}
\email{borya@mccme.ru}

\begin{document}
\maketitle
 \sloppy
\def\pp#1#2{\nfrac{\partial#1}{\partial#2}}
\section*{Abstract}

We integrate the Lifting cocycles
$\Psi_{2n+1},\Psi_{2n+3},\Psi_{2n+5},\dots$ ([Sh1], [Sh2]) on the
Lie algebra $\Dif_n$ of holomorphic differential operators on an
$n$-dimensional complex vector space to the cocycles on the Lie
algebra of holomorphic differential operators on a holomorphic
line bundle~$\lambda$ on an $n$-dimensional complex manifold~$M$
in the sense of Gelfand--Fuks cohomology [GF] (more precisely, we
integrate the cocycles on the sheaves of the Lie algebras of
finite matrices over the corresponding associative algebras). The
main result is the following explicit form of the Feigin--Tsygan
theorem [FT1]:
$$
H^\ndot_\Lie(\gl^\fin_\infty(\Dif_n);\C)=\wedge^\ndot(\Psi_{2n+1},
\Psi_{2n+3},\Psi_{2n+5},\dots).
$$

\section*{Introduction}
The cocycles $\Psi_{2n+1},\Psi_{2n+3},\Psi_{2n+5},\dots$ on the
Lie algebra $\gl_\infty^\fin(\Dif_n)$ of finite matrices over
(polynomial, holomorphic, formal) differential operators on $\C^n$ ($\Psi_i\in\nobreak
C^i_\Lie(\gl^\fin_\infty(\Dif_n);\C)$),
called \emph{Lifting formulas}, were constructed in the author's works
[Sh1], [Sh2].

In the present paper we study the various aspects of the notion
of integral in Lie algebra cohomology, applied to the Lie algebra $\gl_\infty^\fin(\Dif_n)$,
namely, we integrate the cocycles
$\Psi_{2n+1},\Psi_{2n+3},\Psi_{2n+5},\dots$ on this Lie algebra on
an $n$-dimensional complex manifold~$M$. In this way, we obtain
$1$-, $3$-, $5$-, $\dots$ cocycles on the sheaf of the Lie algebras of
holomorphic differential operators
in any holomorphic line bundle~$\lambda$ over $M$ (in
the sense of Section~3).

\subsection{}
It was proved in [FT1] that the cohomology algebra
$H^\ndot(\gl^\fin_\infty(\Dif_n);\C)$ is the
exterior algebra with the generators in degrees $2n+1, 2n+3,
2n+5,\dots$. This result was proved using the spectral sequence,
connecting the Hochschild homology of an associative algebra and
its cyclic homology (and the result of [T]); in particular, any
explicit formulas for the generators do not follow from this
computation.

We prove that
$$
H^\ndot(\gl^\fin_\infty(\Dif_n);\C)\simeq\wedge^\ndot(\Psi_{2n+1},
\Psi_{2n+3},\Psi_{2n+5},\dots),
$$
where $\Psi_{2n+1},\Psi_{2n+3},\Psi_{2n+5},\dots$ are Lifting cocycles
(Theorem~4.3.6 in the case $n=1$ and Theorem~4.4 in the general
case.) We prove this theorem using the notion of the
\emph{integral} in Lie algebra cohomology, which is due to
I.M.\,Gelfand and D.B.\,Fuks [GF]. First of all, let us recall the
classical construction of the Virasoro $2$-cocycle on the Lie
algebra $\Vect(S^{1})$ of smooth vector fields on the circle.

We choose any (formal) coordinate system in each point $x\in
S^{1}$, smoothly depending on the point~$x$; in this way, we obtain the
map $\iota_{x}\colon \Vect(S^{1})\to W_{1}$, connected with each point
$x\in S^{1}$ ($W_{1}$ is the Lie algebra of formal vector field
on the line $\R^{1}$). Let $\Psi_{3}$ be the $3$-cocycle on the Lie
algebra $W_{1}$ (in fact, $\dim H^{0}(W_{1};\C)=\dim
H^{3}(W_{1};\C)=1$, $\dim H^{i}(W_{1};\C)=0$ when $i\ne 0,3$).

We obtain ``in any point $x\in S^{1}$'' the cocycle
$\Psi_{3}(x)=\iota^*_x\Psi_3\in C^3_\Lie(\Vect(S^1);\C)$.
In fact, the cohomological class of all the
cocycles $\Psi_{3}(x)$ is the same,  it does not depend on the
point $x\in S^{1}$ and on the choice of the (formal) coordinate
system in the point~$x$. It follows from this statement, that there
exists a $1$-form $\Theta_{2}$ on $S^{1}$ with the values in
$C^{2}_{\Lie}(\Vect(S^{1});\C)$ such that
$$
d_\DR \Psi_{3}(x)=\delta_{\Lie}\Theta_{2}.
$$
Then $\int_{S^1}\Theta_{2}$ is a cocycle in
$C^2_\Lie(\Vect(S^1);\C)$. Indeed,
$$
\delta_\Lie \int_{S^1}\Theta_2=\int_{S^1}\delta_{\Lie}\Theta_{2}=
\int_{S^1}d_{DR}\Psi_{3}=0.
$$
In fact, the choice of $\Theta_{2}$ is not unique, but there exists
in a sense the canonical choice (see Subsec.~2.2). The cohomological
class $\left[\int_{S^1}\Theta_{2}\right]$ does not depend on the choice of the
coordinate systems.

Let $M$ be an $n$-dimensional complex manifold, $\lambda$ be a
holomorphic line bundle on $M$. Let $\Dif_{\lambda,M}$ be the sheaf of
the associative algebras of holomorphic differential operators in
$\lambda$, and $\D^\ndot_{M}$ be the Dolbeault complex on~$M$:
$$
\D^\ndot_M=\left\{\,0\to
C^\infty_M\nto{\op}\Omega^{0,1}_M\nto{\op}\Omega^{0,2}_M\nto{\op}\dots\,\right\}.
$$

We consider
$
\Gamma_{M}(\gl^\fin_\infty(\Dif_{\lambda,M}\otimes_{\O_M}\D
^\ndot_M))
$
as a $\DG$ Lie algebra (see Remark 3.1).
When $M=\C^n$, the $\DG$ Lie algebra
$
   \Gamma_{\C^n}(\gl^\fin_\infty(\Dif_{\C^n}\otimes_{\O_{\C^n}}\D
^\ndot_{\C^n}))
$
is quasi-isomorphic to $\Dif_n [0]$ as a $\DG$ Lie algebra.

We define the cohomology of the sheaf of Lie algebras $\Dif_{\lambda,M}$
(with the bracket $[a,b]=(a*b-b*a)$) as cohomology of the $\DG$ Lie
algebra $\Gamma_{M}(\Dif_{\lambda,M}\otimes_{\O_M}\D^\ndot_M)$
and cohomology of the sheaf $\gl_\infty^{\fin}(\Dif_{\lambda,M})$
as cohomology of the $\DG$ Lie algebra
$
\Gamma_{M}(\gl^\fin_\infty(\Dif_{\lambda,M}\otimes_{\O_M}\D
^\ndot_M))
$.
The reason is that complex of sheaves
$\Dif_{\lambda,M}\otimes_{\O_{M}}D^\ndot_M$ is quasi-isomorphic
to the sheaf $\Dif_{\lambda,M}[0]$ (as sheaves of associative or Lie
algebras).

For any cocycle $\Psi \in C^{k}_{\Lie}(\gl_\infty^{\fin}(\Dif_{n});\C)$ and
any
singular cycle $\sigma\in C^{\sing}_{l}(M;\C)$ we define (in
Sec.~2 and Subsec.~3.2) the \emph{integral $\int_{\sigma}\Psi$}, which is a
\emph{cocycle} in
$$
C^{k-l}_\Lie(\Gamma_{M}(\gl^\fin_\infty(\Dif_{\lambda,M})\otimes_{\O_M}\D
^\ndot_M)).
$$
The crucial point is that if $[\int_{\sigma}\Psi]\ne 0$ than $[\Psi]\ne 0$
(here $[\dots]$ stands for
the cohomological class) (Theorem 2.3(1)). Therefore the notion of integral
gives us an effective tool to prove the cohomological nontriviality of
cocycles in $C^\ndot_\Lie(\gl^\fin_\infty(\Dif_n);\C)$.

Let $n=1$, $M=\C P^1$. In this simplest case the sheaf of holomorphic
differential operators (in any holomorphic line bundle $\lambda$) has not
higher cohomology (as a sheaf), and the $DG$ Lie algebra
$
\Gamma_{\CP^1}(\gl^\fin_\infty(\Dif_{\lambda,\CP^1})
\otimes_{\O_{\CP^1}}\D^\ndot_{\CP^1})
$
is quasi-isomorphic to the Lie algebra of global differential operators
$\Gamma(\Dif_{\lambda,\C P^1}) [0]$. We will denote the last Lie algebra by
$\Dif_{\lambda,\CP^1}$. (The situation is the same for generalized flag
varieties, in particular, for projective spaces and flag varieties).

Using the method of [FT1] it is not difficult to show that
$H^\ndot_\Lie(\gl^\fin_\infty(\Dif_{\C P^1});\C)$ is the
exterior algebra with the unique generator in dimension $1$ and two generators
in each dimension $3, 5, 7, \dots$.

We calculate the integral
$\int_{\CP^1}\Psi_{2k+1}\in C^{2k-1}_\Lie (\gl^\fin_\infty(\Dif_{\CP^1});\C)$ and
 prove, that it has nonzero value on
the $(2k-1)$-cycle in $H_{2k-1}(\gl^\fin_\infty(\Dif_{\CP^1});\C)$ which is the
image of some $(2k-1)$-cycle in $H_{2k-1}(\gl^\fin_\infty;\C)$
under the inclusion
$\gl^\fin_\infty\hookrightarrow\gl^\fin_\infty(\Dif_{\CP^1})$
($A\mapsto A\otimes 1$) ($k\ge 1$)
(Subsection 4.1--4.3). Then it follows from Theorem 2.3(1) that $\Psi_{2k+1}$
is not cohomologous to zero in
$C^\ndot_\Lie(\gl^\fin_\infty(\Dif_1;\C)$ for any $k\ge 1$. It follows from
this result, calculation of [FT1], and the Hopf algebra structure on
cohomology
$H^\ndot_\Lie(\gl^\fin_\infty(\Dif_1);\C)$ that
$$
H^\ndot_\Lie(\gl^\fin_\infty(\Dif_1);\C)=\wedge^\ndot(\Psi_3,\Psi_
5, \Psi_7,\dots)
$$
(in fact, the cocycles $\Psi_{2k+1}$, $k\ge 1$, are
\emph{primitive} elements with respect to the Hopf algebra structure on
$H^\ndot_\Lie(\gl^\fin_\infty(\Dif_1);\C)$).

Unfortunately, we have not found any simple calculation of the integrals using the
Dolbeault complex as we described above. Our calculation is based on the
\v Cech resolution, and we define the notion of integral from viewpoint of the
\v Cech resolution (in the case of $\C P^n$) in Subsec. 4.1--4.2. We believe,
that
both definitions of the integral coincide.

Let $\iota\colon\Dif_{\CP^1}\hookrightarrow\Dif_1$ be the inclusion defined by
a
choice of coordinate system in a point of $\C P^1$. We prove that
$$
H^\ndot_\Lie(\gl_\infty^\fin(\Dif_{\CP^1});\C)\simeq
\wedge^\ndot\left(\int_{\CP^1}\Psi_3;\ \iota^*\Psi_3,
\int_{\CP^1}\Psi_5;\ \iota^*\Psi_5,\int_{\CP^1}\Psi_7;\dots\ \right).
$$

Using the integration on $\CP^n$, we prove that
$$
H^\ndot_\Lie(\gl_\infty^\fin(\Dif_n);\C)\simeq
\wedge^\ndot(\Psi_{2n+1},\Psi_{2n+3},\Psi_{2n+5},\dots)
$$
for any $n\ge1$. The situation here is quite complicated,
because we don't know any way to calculate the integrals directly.
Let $\lambda=\O(\nl)$, $\nl\in\Z$, and we consider
$$
\int_{\CP^n}\Psi_{2k+1}\in H_\Lie^{2k-1}(\gl_\infty^\fin
(\Dif_{\lambda,\CP^n});\C),\qquad k\ge n.
$$
For a matrix $(2k-2n+1)$-cycle
$
\gamma\in H_{2k-2n+1}(\gl_\infty^\fin;\C)$ we consider
$\int_{\CP^n}\Psi_{2k+1}$ as a polynomial function on~$\nl$.
In fact, this is a polynomial of $n$-th degree. It is
not easy to calculate this polynomial, but we calculate
its leading coefficient and prove, that
it is not equal to zero for some~$\gamma$. Therefore,
$$
\int_{\CP^n}\Psi_{2k+1}\in C_\Lie^{2k-2n+1}
(\gl_\infty^\fin(\Dif_{\lambda,\CP^n});\C)
$$
is not cohomologous to zero, and it follows from the
Theorem~2.3.(1) that
$$
\Psi_{2k+1}\in C_\Lie^{2k+1}(\gl_\infty^\fin(\Dif_n);\C)
$$
is not cohomologous to zero for any $k\ge n$. The
remaining part of the proof is the same as in the case $n=1$.

Finally, let us note that the question on the cohomological
nontriviality of the pull-back $j^*\Psi_{2k+1}$ with respect
to the inclusion $j\colon\Dif_n\hookrightarrow\gl_\infty^\fin(\Dif_n)$
($\D\mapsto E_{11}\otimes\D$) remains open; it solved only in the simplest
case $k=n$ ([Sh1], Sect.~2). We don't know how to check the
nontriviality of these cocycles because it is not known any explicit
formulas for cycles in~$C_\ndot(\Dif_n;\C)$.
So the case of the Lie algebra $\gl_\infty^\fin(\Dif_n)$
turns out more simpler than the case of the Lie algebra~$\Dif_n$
itself, because there exist matrix cycles in the case of
$\gl_\infty^\fin(\Dif_n)$.

On the other hand, the values of Lifting formulas itself
(not their integrals) on the matrix cycles are equal to~$0$,
as well as the values of the pull-back $\iota^*\Psi_{2k+1}$
on the matrix cycles.

\subsection{Content of the paper}

\begin{enum}
\item[Section 1] we recall ([Sh1], [Sh2]) 
the construction of the Lifting cocycles
$\Psi_{2n+1},\Psi_{2n+3},\Psi_{2n+5},\dots$ 
on the Lie algebra $\gl_\infty^\fin(\Dif_n)$;
\item[Section 2] we recall the definition of the notion of integral
in Lie algebra  cohomology in the classical case of the Lie
algebra of smooth vector fields on a real manifold;
all this material is basically standard ([GF]),
the minor modifications are due to Boris Feigin;
\item[Section 3] we give the definition of the integral in the case
of the Lie algebra of holomorphic differential operators
on a holomorphic line bundle~$\lambda$ on a complex manifold~$M$
(more precisely, the sheaf of the Lie algebras, $\dots$);
when $M$ is compact, the holomorphic noncommutative residue
is defined. In this case, the value of the noncommutative
residue on the identity differential operator is equal, up
to a constant depending only on~$M$, to the Euler characteristic
of the bundle~$\lambda$ (Conjecture~3.3);
\item[Section 4] we calculate the integrals in the case $M=\CP^n$,
$n\ge1$, $\lambda=\O(\nl)$. We use the \v Cech resolution
instead of the Dolbeault resolution of Section~3. We prove
here, that
$
H^\ndot_\Lie(\gl_\infty^\fin(\Dif_n);\C)=\wedge^\ndot(\Psi_{2n+1},
\Psi_{2n+3},\Psi_{2n+5},\dots)
$.
\end{enum}

\subsection{Acknowledgements}

I am very much indebted to Boris Feigin, almost
all I know about Lie algebra cohomology
I learned from the numerous informal
discussions with him.

The first ideas of this paper appeared in February 1998,
when I worked in the Institut des Hautes \'Etudes Scientifiques
(Bures--sur--Yvette). I am grateful to the I.H.E.S. for the
hospitality and very stimulating atmosphere.

\section{Lifting formulas ([Sh1], [Sh2])}

\subsection{}

It will be convenient to give the definitions in a bit more generality.

Let $\A$ be an associative algebra, we will also consider
$\A$ as the Lie algebra  with the bracket
$[a,b]=a\cdot b-b\cdot a$. Let $\Tr\colon\A\to\C$ be
a trace on the associative algebra~$\A$ (i.\,e.\
$\Tr[a,b]=0$ for any $a,b\in\A$); we denote by $\Der_\Tr\A$
the Lie algebra of all the derivations~$D$ of the associative
algebra~$\A$ such that $\Tr(Da)=0$ for any
$a\in\A$.

Lifting formulas are the formulas for $(k+1)$-,
$(k+3)$-, $(k+5)$- $\dots$ cocycles on the \emph{Lie algebra~$\A$}
constructed by derivations $D_1,\dots,D_k\in\Der_\Tr\A$
such that the following conditions (i)--(ii) hold:

(i) $[D_i,D_j]=\ad Q_{ij}$\quad($Q_{ij}\in\A$, $Q_{ji}=-Q_{ij}$);

(ii) $\Altl_{i,j,l}D_l(Q_{ij})=0$

\noindent
(Let us note that in any case $\Altl_{i,j,l}D_l(Q_{ij})$
lies in the center of the Lie algebra~$\A$.).

In the simplest case, when $Q_{ij}=0$ for all $i,j$
the formula for $(k+1)$-cocycle is given by the following formula:
\begin{equation}
\Psi_{k+1}(A_1,\dots,A_{k+1})=\Altl_{A_1,\dots,A_{k+1}}
\Altl_{D_1,\dots,D_k}\Tr(D_1A_1\cdot\ldots\cdot D_kA_k\cdot A_{k+1}).
\end{equation}
In the general case, the r.h.s.\ of~(1) is the ``leading term''
of the $(k+1)$-cocycle, also exist terms linear, quadratic, $\dots$,
$\left[\nfrac k2\right]$-th degree in~$Q_{ij}$.

Let $\A=\PDif_n$ be the associative algebra of the formal
pseudo-differential operators on~$(S^1)^n$ ([A]). One can define
$2n$ (exterior) derivations of this algebra: $\ad(\ln x_1)$,
$\dots$, $\ad(\ln x_n)$; $\ad(\partial_1)$, $\dots$, $\ad(\ln\partial_n)$.
We use the formal symmetry between the symbols $\partial_i$
and $x_i$ arising from the generating relation
$[\partial_i,x_i]=1$ to define $\ad(\ln\partial_i)$.
These derivations as well as the corresponding $2$-cocycles
firstly appear in~[KK].
The trace on $\PDif_n$ is the ``noncommutative residue''
([A]): it is defined as the coefficient at the term
$x_1^{-1}\dots x_n^{-1}\partial^{-1}_1\dots\partial_n^{-1}$
in any coordinate system; it was proved in~[A], that
$\Tr[a,b]=0$ for any $a,b\in\A$. It is easy to see
([Sh1]) that in fact $\ad(\ln x_i),\ad(\ln\partial_i)\in
\Der_\Tr\PDif_n$. It was proved in [Sh1] that
\begin{multline}
[\ad(\ln\partial_i),\ad(\ln\partial_i)]=\\
\ad\left(x_i^{-1}\partial_i^{-1}+\nfrac12x_i^{-2}\partial_i^{-2}+
\nfrac23x_i^{-3}\partial^{-3}_i+\ldots+\nfrac{(n-1)!}n
x_i^{-n}\partial_i^{-n}+\dots\right).
\end{multline}
Therefore, the set of $2n$ derivations
$$
\{\ad(\ln x_1),\dots,\ad(\ln x_n);\ad(\ln\partial_1),\dots,
\ad(\ln\partial_n)\}
$$
satisfy the conditions (i), (ii) above.

We can reply this construction in the case $\A=\gl_\infty^\fin(\PDif_n)$,
here
$$
\Tr_{\gl_\infty^\fin(\PDif_n)}=\Tr_{\PDif_n}\circ\Tr_{\gl_\infty^\fin}.
$$
The derivations $\ad(\ln x_1)$, $\dots$, $\ad(\ln x_n)$;
$\ad(\ln\partial_1)$, $\dots$, $\ad(\ln\partial_n)$
act in the obvious way on the algebra
$\gl_\infty^\fin(\PDif_n)$, and the conditions (i)--(ii)
hold again.

\begin{remark}
Strictly speaking, $[\ad(\ln\partial_i),\ad(\ln x_i)]=\ad Q$,
where $Q$ is an infinite matrix in the case of the Lie algebra 
$\gl_\infty^\fin(\Dif_n)$
.
\end{remark}

By definition, the Lifting formulas for the algebra $\Dif_n$
are pull-backs of the Lifting formulas for $\PDif_n$
with respect to the natural imbedding $\Dif_n\hookrightarrow
\PDif_n$, and also for the Lifting formulas in the case of the algebra
$\gl_\infty^\fin(\Dif_n)$.

We construct the Lifting formulas for $\gl_\infty^\fin(\Dif_1)$
($k=2$) in Subsection~1.2 and for general~$k$
in Subsection~1.3.

\subsection{The case $n=1$ {\rm([Sh1], [KLR])}}

Let $k=2$, so we have two derivations~$D_1,D_2$ such that
$[D_1,D_2]=\ad Q$.

The simplest Lifting formula for $3$-cocycle is given by the
following formula:
\begin{equation}
\Psi_3(A_1,A_2,A_3)=\Altl_{A,D}\Tr(D_1A_1\cdot D_2A_2\cdot A_3)
+\Altl_A\Tr(Q\cdot A_1\cdot A_2\cdot A_3)
\end{equation}

In the general case, for $i=2,3,4,\dots$
\begin{multline}
\Psi_{2i+1}(A_1,\dots,A_{2i+1})=\Altl_A\Tr(Q\cdot A_1\cdot A_2
\cdot\ldots\cdot A_{2i+1})+\\
+\Altl_{A,D}\Tr(D_1A_1\cdot D_2A_2\cdot A_3\cdot\ldots\cdot A_{2i+1}+
D_1A_1\cdot A_2\cdot A_3\cdot D_2A_4\cdot A_5\cdot\ldots\cdot A_{2i+1}+\\
+D_1A_1\cdot A_2\cdot A_3\cdot A_4\cdot A_5\cdot D_2A_6\cdot A_7
\cdot\ldots\cdot A_{2i+1}+\ldots+\\
+a_s\cdot D_1A_1\cdot A_2\cdot\ldots\cdot A_s\cdot
D_2A_{s+1}\cdot A_{s+2}\cdot\ldots\cdot A_{2i+1})
\end{multline}
$$
\text{where}\qquad
\begin{cases}
s=i+1,\quad a_s=\nfrac12\qquad&\text{if $i$ is even}\\
s=i,\quad a_s=1\qquad&\text{if $i$ is odd}
\end{cases}\ .
$$

\begin{theorem}[{[Sh1]}, Sect.~1, {[KLR]}]
The formulas for $\Psi_3,\Psi_5,\Psi_7,\dots$
given by \emph{(3), (4)} define cocycles on the Lie algebra~$\A$.
\qed
\end{theorem}

\subsection{The general case}

In this Subsection we suppose that
$D_1,\dots,D_k\in\Der_\Tr\A$ satisfy conditions (i), (ii)
from Subsec.~1.1.

First of all, let us write the formula for $5$-cocycle in the case
$k=4$ (this is the simplest Lifting formula in this case):
\begin{multline}
\Psi_5(A_1,A_2,A_3,A_4,A_5)=\\
\begin{aligned}\relax
&=\Altl_{A,D}\Tr\{D_1A_1\cdot D_2A_2\cdot D_3A_3\cdot D_4A_4
\cdot A_5&&\\
&
\begin{aligned}\relax
&+A_1\cdot Q_{12}\cdot A_2\cdot D_3A_3\cdot D_4A_4\cdot A_5\\
&+D_1A_1\cdot A_2\cdot Q_{23}\cdot A_3\cdot D_4A_4\cdot A_5\\
&+D_1A_1\cdot D_2A_2\cdot A_3\cdot Q_{34}\cdot A_4\cdot A_5
\end{aligned}&
\left.\vphantom{
\begin{aligned}\relax
&+A_1\cdot Q_{12}\cdot A_2\cdot D_3A_3\cdot D_4A_4\cdot A_5\\
&+D_1A_1\cdot A_2\cdot Q_{23}\cdot A_3\cdot D_4A_4\cdot A_5\\
&+D_1A_1\cdot D_2A_2\cdot A_3\cdot Q_{34}\cdot A_4\cdot A_5
\end{aligned}}
\right]&\quad\text{terms, linear in $Q_{ij}$}\\
&\,+A_1\cdot Q_{12}\cdot A_2\cdot A_3\cdot Q_{34}\cdot A_4\cdot
A_5\}&]\,&\quad\text{term, quadratic in $Q_{ij}$}
\end{aligned}
\end{multline}
We don't alternate the symbols $i,j$ in $Q_{ij}$.

The next step is the formula for $(k+1)$-cocycle for any~$k$.
We need some preparations. We consider the interval of the
length $k-2$ with some \emph{marked} points, such that the
\emph{distance between any two marked points is greater or
equal than~$2$}.
Let us denote by $I_l$ the set of all such marked intervals
with $l$ marked points ($1\le l\le\left[\nfrac k2\right]$).
Denote by $1,\dots,k-1$ the integral points of the interval.

\begin{defin}
Suppose that $t\in I_l$ and $i_1<\ldots<i_l$
are its marked points ($1\le i_1$,\ \ $i_k\le\left[\nfrac k2\right]$
and $i_{s+1}-i_s\ge2$ for all $s=1,\dots,l-1$). Then
$$
\O(t)=\Altl_{A,D}\Tr(P_{1,t}\cdot\ldots\cdot P_{k+1,t})
$$
where
\begin{align*}\relax
&\left.
\begin{aligned}\relax
&P_{j,t}=A_j\cdot Q_{j,j+1}\\
&P_{j+1,t}=A_{j+1}
\end{aligned}\right\}\quad&\parbox{6.5cm}{when $j=i_s$ for
some $s=1,\dots,l$, i.\,e., if point $j$ is marked}\\
&\left.P_{j,t}=D_jA_j\right.\quad&\parbox{6.5cm}{when $j$ and $j-1$ are not
marked and $j\ne k+1$}\\
&\left.P_{k+1,t}=A_{k+1}\right.&
\end{align*}
\end{defin}

\begin{example}
If $k=6$ and $t=\vcenter{\epsfbox{pic.7}}$, $t\in I_2$ then
$$
\O(t)=\Altl_{A,D}\Tr(A_1\cdot Q_{12}\cdot A_2\cdot D_3A_3
\cdot D_4A_4\cdot D_5A_5\cdot D_6A_6\cdot A_7).
$$
If $t=\vcenter{\epsfbox{pict.8}}$, $t\in I_3$,
then
$$
\O(t)=\Altl_{A,D}\Tr(A_1\cdot Q_{12}\cdot A_2\cdot A_3\cdot Q_{34}
\cdot A_4\cdot A_5\cdot Q_{56}\cdot A_6\cdot A_7).
$$
\emph{We don't alternate the symbols $i,j$ in $Q_{ij}$
in this definition}.
\end{example}

\begin{theorem}
Let $\Sigma_l=\sum\limits_{t\in I_l}\O(t)$. Then
\begin{equation}
\Psi_{k+1}=\Altl_{A,D}\Tr(D_1A_1\cdot D_2A_2\cdot
\ldots\cdot D_kA_k\cdot A_{k+1})+\Sigma_1+
\Sigma_2+\ldots+\Sigma_{\left[\frac k2\right]}
\end{equation}
is a $(k+1)$-cocycle on the Lie algebra~$\A$.
\end{theorem}
\qed

This Theorem was proven in [Sh2], Section~3.

It remains to define the cocycles $\Psi_{k+3},\Psi_{k+5},\Psi_{k+7},
\dots$ (for $k$ derivations $D_1,\dots,D_k$).
First of all, let us suppose that $Q_{ij}=0$ for all $i,j$.

\begin{defin}
The set $a^s_\even$ is the set of the sequences $\{a_i\}$ of
the length $k+2s$ such that:

(i) $a_i\in\{0,1\}$ for all $i=1,\dots,k+2s$;

(ii) $a_1=1$;

(iii) $k$ of the $a_i$'s are equal to~$1$ and $2s$ of the $a_i$'s
are equal to~$0$.

(iv) the number of $0$'s between the two nearest $1$'s is
\emph{even}; this condition also should hold for the ``tail''
of the sequence~$\{a_i\}$, as if the number's~$a_i$ were
placed on a circle.
\end{defin}

We define the $(k+2s-1)$-cochain
$R_{a_1,\dots,a_{k+2s}}(A_1,\dots,A_{k+2s-1})$ for each
$\{a_i\}\in a_\even$. Roughly speaking, we just shorten any
sequence of consecutive zeros to a single zero; after this
procedure the sequence will have odd length. We choose the
first such sequence of the consecutive $0$'s.
More precisely, let us define the sequence
$\{\wa_i\}\in a_\even$ in the following way:

let $s_1=\min\limits_i(a_i=0)$, $s_2=\max\limits_{i>s_1}(a_i=1)$; then:
\begin{align*}\relax
&\wa_1,\dots,\wa_{s_1-1}=1;\\
&\wa_{s_1},\dots,\wa_{s_2-2}=0;\\
&\wa_i=a_{i+1}\quad\text{for}\quad s_2-1\le i\le k+2s-1.
\end{align*}

\begin{defin}
$$
R_{a_1,\dots,a_{k+2s}}(A_1,\dots,A_{k+2s-1})=
\Altl_{A,D}\Tr(P_1\cdot\ldots\cdot P_{k+2s-1}),
$$
where
\begin{align*}\relax
&P_i=D_{j(i)}A_i\quad\text{for}\quad\wa_i=1,\\
&P_i=A_i\quad\text{for}\quad\wa_i=0
\end{align*}
and $j(i)$ is defined in the following way: $j(1)=1$ and
$j(i_1)<j(i_2)$ when $i_1<i_2$ and $j=1,\dots,n$. In other words,
$j$~takes values from $1$ to~$k$ \emph{in turn}.
\end{defin}

\begin{theorem}[{[Sh2]}, Section 1]
\begin{equation}
\Psi^0_{k+2s-1}=\sum_{\{a_i\}\in
a_\even}(-1)^{s_1}R_{a_1,\dots,a_{k+2s}}(A_1,\dots,A_{k+2s-1})
\end{equation}
is a $(k+2s-1)$-cocycle on the Lie algebra~$\A$ when all $Q_{ij}=0$.
\end{theorem}
\qed

The formula (7) is the generalization of the formula (1) on the
case of the higher cocycles.

Let us consider the general case, $Q_{ij}\ne0$. We quantize each
summand $R_{a_1,\dots,a_{k+2s}}$ separately.

\begin{defin}
(i) Denote by $\Circle_l^{a_1,\dots,a_{k+2s}}$
the set of all circles with $k+2s-1$ integral points from
which~$l$ ($1\le l\le\left[\nfrac k2\right]$) are marked. The
distance between any two marked points is $\ge2$. The points are enumerated by
$1,\dots,k+2s-1$. Point $i$ may be marked only if $\wa_i=1$ and
$\wa_{i+1}=1$ (we suppose that $\wa_{k+2s}=\wa_1$).

(ii) for $t\in\Circle_l^{a_1,\dots,a_{k+2s}}$
we define $\O(t)$ by the analogy with the definition above after
replacing the interval with the circle. We don't alternate the
symbols $i,j$ in $Q_{ij}$ in this definition.

(iii)
$$
\Sigma_l^{a_1,\dots,a_{k+2s}}=\sum_{t\in\Circle_l^{a_1,\dots,a_{k+2s}}}\O(t)
$$
\end{defin}

\begin{theorem}[{[Sh2]}, Section 3]
The formula
\begin{equation}
\Psi_{k+2s-1}=\Psi_{k+2s-1}^0+\sum_{\{a_i\}\in a_\even}(-1)^{s_1}
\sum_{l\ge1}\Sigma_l^{a_1,\dots,a_{k+2s}}
\end{equation}
defines a $(k+2s-1)$-cocycle on the Lie algebra~$\A$, where
$\Psi_{k+2s-1}^0$ is defined by the formula~\emph{(7)}.
\end{theorem}
\qed

\section{Integration in the Lie algebra cohomology [GF]}

\subsection{$C^\infty$-case}

Let $M$ be an $n$-dimensional $C^\infty$-manifold, let $W_n$ be
the Lie algebra of $C^\infty$-vector fields on~$\R^n$ and let
$\Vect(M)$ be the Lie algebra of $C^\infty$-vector fields on~$M$.

The integration is a procedure, corresponding to any cocycle
$\xi\in C_\Lie^k(W_n;\C)$ and any singular cycle $\sigma\in
C^\sing_l(M,\C)$ the cocycle $\int_\sigma\xi\in
C_\Lie^{k-l}(\Vect(M);\C)$ such that:

(i) $\left[\int_\sigma\xi\right]=0$\quad if $[\xi]=0$;

(ii) $\left[\int_\sigma\xi\right]=0$\quad if $[\sigma]=0$

\noindent
(here $[\dots]$ denotes the cohomological class).

Let $x\in M$ be any point of $M$, and let $\phi\colon U\to M$ be
any coordinate system in the point~$x$\ \ ($U$ is a neighbourhood
of~$0$ in~$\R^n$ and $\phi(0)=x$). The map $\phi$ induces the map
$\phi_\Vect\colon\Vect(M)\to W_n$, and the corresponding map of
the cochain complexes:
$$
\phi^*_\Vect\colon C_\Lie^\ndot(W_n;\C)\to
C_\Lie^\ndot(\Vect(M);\C).
$$

\begin{lemma}
Let $\xi\in C_\Lie^\ndot(W_n;\C)$ be a cocycle. Then the
cohomological class $[\phi^*_\Vect(\xi)]$ does not depend on the
choice of the point $x\in M$ and of the coordinate system
$\phi\colon U\to M$.
\end{lemma}

\begin{proof}
The $\ad$-action of the Lie algebra $\Vect(M)$ moves points
of~$M$ and induces the trivial action in cohomology.
\end{proof}

Let $\xi\in C^k_\Lie(W_n;\C)$ be a cocycle. We choose a
coordinate system $\phi_x\colon U\to M$\ \ ($\phi_x(0)=x$) in all
the points $x\in M$, \emph{smoothly depending on the
point~$x$}; by definition,
$$
\xi(x)=\phi_{x,\Vect}^*(\xi)\in C_\Lie^k(\Vect(M);\C).
$$
According to Lemma, all the cocycles~$\xi(x)$ ($x\in M$) are
cohomologous to each other. Therefore, there exists an element
$\Theta_1\in\Omega_M^1\otimes C_\Lie^{k-1}(\Vect(M);\C)$
such that
$d_{\DR}\xi(x)=\delta_\Lie\Theta_1$ (here $d_\DR$ denotes the de
Rham differential and $\delta_\Lie$ denotes the differential in
the cochain complex).

Furthermore, we will find elements $\Theta_i\in\Omega^i_M\otimes
C_\Lie^{k-i}(\Vect(M);\C)$ such that
\begin{equation}
\begin{cases}
d_\DR\xi(x)=\delta_\Lie\Theta_1\\
d_\DR\Theta_1=\delta_\Lie\Theta_2\\
\hbox to 3cm{\dotfill}\\
d_\DR\Theta_{n-1}=\delta_\Lie\Theta_n
\end{cases}\ .
\end{equation}
It is obvious that if solution of~(9) exist, it is not unique.

\begin{propos}
For any  singular cycle $\sigma\in C_l^\sing(M;\C)$ and for
any solution $\{\Theta_i\}$ of the system \emph{(9)} the integral
of $l$-form
$$
\int_\sigma\Theta_l\in C_\Lie^{k-l}(\Vect(M);\C)
$$
is a cocycle.
\end{propos}

\begin{proof}
$$
\delta_\Lie\int_\sigma\Theta_l=\int_\sigma\delta_\Lie\Theta_l=
\int_\sigma d_\DR\Theta_{l-1}\overset{\text{by
Stokes formula}}{=}0.\qed
$$
\let\qed\relax
\end{proof}

It is obvious (see Subsec.~2.3) that the condition~(ii) holds for
any solution of the system~(9). We have in mind (but we don't
make the exact statement) that solution of~(9) is unique if we
require the execution of the condition~(i).

\subsection{The solution of (9)}

In this Subsection we construct the canonical solution of the
system~(9).

Let $x\in M$ and we have choose a coordinate system in any
point, smoothly depending on the point. Then any tangent vector
$v\in T_xM$ in the point~$x$ determines an infinitesimal
$1$-parametric group of diffeomorphisms on~$\R^n$, i.\,e. the
element $t_v\in W_n$.

In the cochain complex of any Lie algebra~$\g$ the following
identity holds:
\begin{equation}
\delta_\Lie\circ\iota_t\pm\iota_t\circ\delta_\Lie=\ad(t)
\end{equation}
($t\in \g$ and $\iota_t\colon C^\ndot_\Lie(\g;\C)\to
C^{\ndot-1}_\Lie(\g;\C)$ is the substitution of $t$ for the first
argument).

Let $\widetilde t_v$ be any $C^\infty$-vector field on $M$ such that
$\widetilde t_v|_\varphi (u)=\varphi_{x,\Vect}(t_v)$. It follows
directly from (10) that
\begin{equation}
\delta_\Lie\circ\iota_{\widetilde t_v}\pm \iota_{\widetilde t_v}
\circ\delta_\Lie=\ad (\widetilde t_v)
\end{equation}
in the cochain complex $C^\ndot_\Lie(\Vect(M);\C).$

\begin{lemma}
\begin{equation}
(d_{\DR}\xi(x))v=\ad\widetilde t_v(\xi(x)).
\end{equation}
\end{lemma}

\begin{remark}
The r.h.s. of (12) does not depend on prolongation
$\widetilde t_v$ of vector field~$t_v$.
\end{remark}

\begin{proof}
The derivation in the direction of the vector $v$ corresponds to
the infinitesimal $1$-parametric group of diffeomorphisms
and $\ad(\widetilde t_v)$ corresponds to the $\Ad$-action of the
diffeomorphism.
\end{proof}

It follows from (11), (12) that
$$
(d_{\DR}\xi(x))v=\delta_\Lie\circ\iota_{\widetilde t_v} \xi(x)
$$
or, in other words, if we set
\begin{equation}
\Theta_1 (v)(A_1,\dots,A_{k-1})=\xi(x)(\widetilde t_v,
A_1,\dots,A_{k-1})
\end{equation}
then $(d_{\DR}\xi(x))v=\delta_\Lie \Theta_1(v)$.

It is obvious that $\Theta_1$ does not depend on the
prolongation $\widetilde t_v$ of the vector field $t_v$.

Furthermore, let us set
\begin{equation}
\Theta_i(v_1,\dots,v_i)(A_1,\dots,A_{k-i})=\xi(x)
(\widetilde t_{v_1},\dots,\widetilde t_{v_i},A_1,\dots,A_{k-i}).
\end{equation}

\begin{theorem}
The elements $\{\Theta_i\}$ defined by the formula~\emph{(14)},
$\Theta_i\in\Omega_M^i\otimes C_\Lie^{k-i}(\Vect(M);\C)$, gives
us a solution of the system~\emph{(9)}.
\end{theorem}

\begin{proof}
We have proved the statement in the case $i=1$. We prove in the
case $i=2$, the general case in analogous.

By the Cartan formula, we have:
\begin{equation}
d_{\DR}(\Theta_1)(v_1(x),v_2(x))=\Theta_1([v_1,v_2](x))
-v_1(\Theta_1(v_2(x))+v_2(\Theta_1(v_1(x)).
\end{equation}

In the formula (15) $v_1$ and $v_2$ are vector fields such that
$v_i|_x=v_i(x)$.

According to Lemma, we have:
\begin{gather}
v_1(\Theta_1(v_2))(x)=\ad(\widetilde t_{v_1})(\Theta_1(v_2(x)),\\
v_2(\Theta_1(v_1))(x)=\ad(\widetilde t_{v_2})(\Theta_1(v_1(x))
\end{gather}
if we choose the vector fields $v_1$ and $v_2$ such that
$[v_1,v_2]=0$. Note also that $[\widetilde t_{v_1(x)},\widetilde
t_{v_2(x)}]=0$ for any two vectors $v_1(x),v_2(x)$. Therefore
\begin{equation}
d_{\DR}(\Theta_1)(v_1,v_2)=-[\widetilde t_{v_1},\xi(x)(\widetilde
t_{v_2},A_1,\dots,A_{k-1})]+[\widetilde
t_{v_2},\xi(x)(t_{v_1},A_1,\dots,A_{k-1})].
\end{equation}
Now it follows from the cocycle condition for $\xi(x)$ that the
r.h.s.\ of~(18) is equal to $\delta_\Lie\Theta_2(v_1,v_2)$.
\end{proof}

\subsection{Conditions (i) and (ii)}

\begin{theorem}
\emph{(1)} If $[\xi]=0$ then $\left[\int_\sigma\xi\right]=0$ for
any manifold~$M$ and $\sigma\in C_\ndot^\sing(M;\C)$.

\emph{(2)} If $[\sigma]=0$ then $\left[\int_\sigma\xi\right]=0$
for any~$\xi$.
\end{theorem}

\begin{proof}
(1): Let $\xi=\delta_\Lie\eta$. Then
$$
\Theta_i(v_1,\dots,v_i)(A_1,\dots,A_{k-i})=(\delta_\Lie\eta)(\widetilde
t_{v_1},\dots,\widetilde t_{v_i}, A_1,\dots,A_{k-i}).
$$
We have:
$$
(\delta_\Lie\eta)(\widetilde t_{v_1},\dots,\widetilde
t_{v_i},A_1,\dots,A_{k-i})=\delta_\Lie(\eta(\widetilde
t_{v_1},\dots,\widetilde t_{v_i},A_1,\dots,A_{k-i})
$$
modulo the following expressions:

(i) $\ad(\widetilde t_{v_j})\eta(\widetilde t_{v_1},\dots,\widehat{\widetilde
t}_{v_j},\dots,\widetilde t_{v_i},A_1,\dots,A_{k-i})$;

(ii) $\eta([\widetilde t_{v_\alpha},\widetilde t_{v_\beta}],\widetilde
t_{v_1},\dots,\widehat{\widetilde t}_{v_\alpha},\dots,\widehat{\widetilde
t}_{v_\beta},\dots,\widetilde t_{v_i},A_1,\dots,A_{k-i})$.

By formula~(12), the summands of the type~(i) are equal to
$(d_{\DR}\eta)(v_j)$ and therefore their pull-backs on the
singular cycle~$\sigma$ are equal to~$0$.

On the other hand, $[\widetilde t_{v_\alpha},\widetilde t_{v_\beta}]=0$
for any $v_\alpha$ and~$v_\beta$, because the corresponding
\hbox{$1$-parametric} groups commute with each other.

Therefore, the summands of both types (i) and (ii) are equal
to~$0$ after pull-backs on~$\sigma$, and
$$
\int_\sigma\Theta_i=\int_\sigma\delta_\Lie(\eta(\widetilde
t_{v_1},\dots,\widetilde
t_{v_i},A_1,\dots,A_{k-i}))=\delta_\Lie\int_\sigma\eta(\widetilde
t_{v_1},\dots,\widetilde t_{v_i},A_1,\dots,A_{k-i}).
$$

(2): $\int_{\partial\gamma}\Theta_i=\int_\gamma
d_{\DR}\Theta_i=\int_\gamma\delta_\Lie\Theta_{i+1}=
\delta_\Lie\int_\gamma\Theta_{i+1}$.
\end{proof}

\begin{corollary}
The cohomological class of the integral $\int_\sigma\Theta_i$
does not depend on the
choice of the coordinate systems \emph($\sigma$~is a singular
$i$-cycle in~$M$\emph).
\end{corollary}

\begin{proof}
It follows directly from Theorem 2.3.(1) and Lemma~2.1.
\end{proof}

\section{Integration in the Holomorphic Case}

\subsection{Extended Lifting Formulas}

\subsubsection{} Let us denote by $\D^\ndot_M$ the Dolbeault
complex of the structure sheaf $\O_M$ of a complex manifold $M$,
$$
\D^\ndot_M=\left\{\,0\to
C^\infty_M\nto{\op}\Omega^{0,1}_M\nto{\op}\Omega^{0,2}_M\nto{\op}\dots\,\right
\};
$$
it is clear, that $\D^\ndot_M$ is a complex of sheaves of $\O$-modules.
We consider the complex of global sections
$\Gamma_M(\D^\ndot_M)$ as a super-commutative associative $\DG$
algebra.

Let us denote by $\Dif_{\lambda,M}$ the sheaf of holomorphic
differential operators in a holomorphic line bundle $\lambda$ on
$M$, then the ``Dolbeault complex of the sheaf
$\Dif_{\lambda,M}$'' is defined as
$\Dif_{\lambda,M}\otimes_{\O_M}\D^\ndot_M$. Also the Dolbeault
complex $\gl^\fin_\infty(\Dif_{\lambda,M})\otimes_{\O_M}\D^\ndot_M$
is defined. On the complex of global sections
$\Gamma_M(\gl^\fin_\infty(\Dif_{\lambda,M})\otimes_{\O_M}\D^\ndot_M)$ the
structure of $\DG$ Lie algebra is defined. In the case $M=\C^n$ or
$M=\CP^n$ the $\DG$ Lie algebra
$\Gamma_M(\gl^\fin_\infty(\Dif_{\lambda,M})\otimes_{\O_M}\D^\ndot_M)$ is
quasi-isomorphic to the ($\DG$) Lie algebra
$\gl^\fin_\infty(\Gamma_M(\Dif_{\lambda,M})) [0]$ as a $\DG$ Lie
algebra.

\begin{remark}
Let $\g$ be a Lie algebra, and let $m^\ndot$ be super-commutative
associative $\DG$ algebra. Then on $\g\otimes_\C m^\ndot$ the
structure of $\DG$ Lie algebra is defined:
$$
[g_1\otimes m_1, g_2\otimes m_2]=[g_1,g_2]\otimes m_1 m_2.
$$
If
$\g$ is the Lie algebra constructed from an associative algebra
$\g$ (with the bracket $[a,b]=a*b-b*a$) then except the above
construction, $\g\otimes m^\ndot$ is an associative $\DG$ algebra
with the product
$(g_1\otimes m_1)*(g_2\otimes m_2)=(g_1*g_2)\otimes(m_1m_2)$, and
there arises the structure of purely even Lie algebra on
$$
g\otimes m^\ndot\colon[g_1\otimes m_1, g_2\otimes
m_2]=(g_1*g_2-(-1)^{\deg m_1\cdot\deg m_2} g_2*g_1)\otimes(m_1m_2).
$$
Everywhere in this work we have in mind the \emph{first}
construction, and we consider $\g\otimes m^\ndot$ as \emph{$\DG$
Lie algebra}.
\end{remark}

\subsubsection{} In the case $M=\C^n$ the $\DG$ Lie algebra
$
\Gamma_{\C^n}(\gl^\fin_\infty(\Dif_{\C^n})\otimes_{\O_{\C^n}}\D^\ndot_{\C^n})
$
is quasi-isomorphic to the $\DG$ Lie algebra
$\gl^\fin_\infty(\Dif_n) [0]$, and there exists the canonical
imbedding of $\DG$ Lie algebras:
$$
\iota\colon\gl^\fin_\infty(\Dif_n) [0]\hookrightarrow
\gl^\fin_\infty(\Dif_n)\otimes_{\O_{\C^n}} \D^\ndot_{\C^n}.
$$
We want to find formulas for the ``extended'' Lifting cocycles
$\hPsi_{2k+1}$ $(k\ge n)$ of the $\DG$ Lie algebra
$\gl^\fin_\infty(\Dif_n)\otimes_{\O_{\C^n}}\D^\ndot_{\C^n}$ such
that $\iota^*(\hPsi_{2k+1})=\Psi_{2k+1}$. To do this, let us note
that there exists another map of $\DG$ Lie algebras:
$$
p\colon\gl^\fin_\infty(\Dif_n)\otimes_{\O_{\C^n}}\D^\ndot_{\C^n}\to
\gl^\fin_\infty(\Dif_n) [0].
$$
By the definition, $p=0$ in the degrees $\ne 0$, and
\begin{equation}
p(\D\otimes_\O f(z_1,\dots,z_n,\overline z_1,\dots,\overline
z_n)=D\otimes_\O f(z_1,\dots,z_n,0,\dots,0)=f(z_1,\dots,z_n)\cdot\D.
\end{equation}

The composition $p\circ \iota=\id_{\gl^\fin_\infty(\Dif_n) [0]}$,
and therefore
\begin{equation}
\hPsi_{2k+1}=p^*\Psi_{2k+1}
\end{equation}
satisfies the condition $\iota^*\hPsi_{2k+1}=\Psi_{2k+1}$.

\subsection{The Integration in the Complex Case}

Let $M$ be a complex manifold of dimension $n$ and let
$\Psi_{2n+1},\Psi_{2n+3},\Psi_{2n+5},\dots$ be Lifting cocycles on
the Lie algebra $\gl^\fin_\infty(\Dif_n)$ (see Section~1). We
choose a holomorphic coordinate system
$\varphi_x\colon U\to M$ ($U\subset\C^n$, $\varphi_x(0)=x$) in any point
$x\in M$, depending smoothly on the point $x\in M$. Then the
cocycle
$$
\varphi^*_x\hPsi_{2k+1}\in
C^{2k+1}_\Lie(\Gamma_M(\gl^\fin_\infty(\Dif_M)\otimes_{\O_M}\D^\ndot_M))
$$
is defined in any point $x\in M$. As in Section 2, the
cohomological class $[\varphi^*_x(\hPsi_{2k+1})]$ does not
depend on the point $x\in M$ and on the holomorphic coordinate
system $\varphi_x$ in the point~$x$.

Analogously, let $\lambda$ be a holomorphic line bundle on
$M$, $p\colon \lambda\to M$ be the projection. Choosing a coordinate
system $\varphi_{x,\lambda}\colon U\times\C\simto p^{-1}(\varphi_x(U))$ in any
point $x\in M$, smoothly depending on the point $x$, we define
the cocycle
$$
\varphi^*_{x,\lambda}(\hPsi_{2k+1})\in
C^{2k+1}_\Lie(\Gamma_M(\gl^\fin_\infty(\Dif_{\lambda,M})\otimes_{\O_M}\D^\ndot_
M)).
$$
Again, the cohomological class
$[\varphi^*_{x,\lambda}(\hPsi_{2k+1})]$ does not depend on the
point $x\in M$ and the coordinate system
$\varphi_{x,\lambda}\colon U\times\C\simto p^{-1}(\varphi_x(U)).$

Furthermore, we define the $i$-form $\Theta_i$ on $M$ (as a
$C^\infty$-manifold) with values in
$$
C^{2k+1-i}_\Lie(\Gamma_M(\gl^\fin_\infty(\Dif_{\lambda,M})
\otimes_{\O_M}\D^\ndot_M))
$$
by the formula:
\begin{equation}
\Theta_i(v_1,\dots,v_i)(\D_1,\dots,\D_{2k+1-i})=\varphi^*_{x,\lambda}(\hPsi_
{2k+1})(\widetilde t_{v_1}\otimes \Id,\dots,\widetilde t_{v_i}\otimes \Id,
\D_1,\dots,\D_{2k+1-i})
\end{equation}
where:

(i) $x\in M$;

(ii) $v_1,\dots,v_i$ are tangent vectors to $M$ (as a $C^\infty$-%
manifold) in point $x$;

(iii) $\widetilde t_{v_1},\dots,\widetilde t_{v_i}$ are differential
operators of order $\le 1$ on $M$ defined as in Subsec.~2.2 and
$\widetilde t_{v_1}\otimes \Id,\dots,\widetilde t_{v_i}\otimes \Id$ are
\emph{infinite} matrices; strictly speaking, they do not lie in the algebra
$\gl^\fin_\infty(\Dif_{\lambda,M})\otimes_{\O_M}C^\infty_M$;

(iv)
$$
\D_1,\dots,\D_{2k+1-i}\in \Gamma_M(\gl^\fin_\infty(\Dif_{\lambda,M})\otimes
_{\O_M}C^\infty_M)=[\Gamma_M(\gl^\fin_\infty(\Dif_{\lambda,M})\otimes_{\O_M}
\D^\ndot_M)]^0.
$$

(v) cochain $\Theta_i(v_1,\dots,v_i)$ is equal to zero when one
of the arguments has grading $\ne 0$.

\begin{lemma}
The cocycles $\Theta_1,\dots,\Theta_{2n}$ satisfy the
system~\emph{(9)}.
\end{lemma}

\begin{proof}
The unique new point (after Theorem~2.2) is the existence of the
component $\delta^{\overline\partial}_\Lie$ in the cochain differential
$\delta_\Lie$, connected with the differential $\overline\partial$ in the
$\DG$ Lie algebra
$\Gamma_M(\gl^\fin_\infty(\Dif_{\lambda,M})\otimes_{\O_M}\D^\ndot_M)$.
However, it follows from the definition (20) of the extended
Lifting cocycles $\hPsi_{2k+1}$ that
$\delta^{\overline\partial}_\Lie(\Theta_i)(v_1,\dots,v_i)\equiv 0$.
\end{proof}

The direct consequence of this Lemma is the statement that
$\int_\sigma \Theta_i$ is a \emph{cocycle} in
$$
C^{2k+1-i}_\Lie(\Gamma_M(\gl^\fin_\infty(\Dif_{\lambda,M})
\otimes_{\O_M}\D^\ndot_M))
$$
for any singular $i$-cycle $\sigma$:
$$
\delta_\Lie\int_\sigma\Theta_i=\int_\sigma\delta_\Lie\Theta_i=
\int_\sigma d_\DR\Theta_{i-1}=0.
$$

One can prove also the direct analogs of Theorem~2.3 and
Corollary~2.3.

\subsection{Holomorphic Noncommutative Residue}

Let $M$ be a \emph{compact} complex manifold, $\dim M=n$,
$\langle M\rangle$
be the fundamental class of $M$, $\langle M\rangle\in H^\sing_{2n}(M;C)$. Then,
according to previous Subsection,
$$
\int_{\langle M\rangle}\Theta_{2n}\in
C^1_\Lie(\Gamma_M(\gl^\fin_\infty(\Dif_{\lambda,M})\otimes_{\O_M}\D^\ndot_M))
$$
is a \emph{cocycle}.

Let us denote
$$
A^\ndot=\Gamma_M(\gl^\fin_\infty(\Dif_{\lambda,M})\otimes_{\O_M}\D^\ndot_M).
$$

\begin{lemma}
Cocycle $\int_{\langle M\rangle}\Theta_{2n}$ defines a linear functional on
$$
\Ker\{\,\overline\partial\colon A_0/[A_0,A_0]\to A_1/[A_0,A_1]\,\}.
$$
\end{lemma}

\begin{proof}
It follows directly from the generalization of the Theorem from
[T] which states that
$$
H^\ndot_\Lie(\gl^\fin_\infty(A);C)\simeq
S^\ndot_\super((HC_\ndot(A)[1])^*)
$$
in the case of a $\DG$ associative algebra $A$.
\end{proof}

When $\D\in\Ker\{\,\overline\partial\colon A_0/[A_0,A_0]\to A_1/[A_0,A_1]\,\}$ we
denote by $\Tr_\lambda(\D)$ this noncommutative residue; in
particular, $\Tr_\lambda(\boldsymbol1_\lambda)$ is well-defined, where
$\boldsymbol1_\lambda$ is the identity differential operator in the line
bundle $\lambda$.

We have proved (Corollary 2.3) that the number
$\Tr_\lambda(\boldsymbol1_\lambda)$ does not depend on the choice of the
coordinate systems in the definition of $\int_{\langle M\rangle}\Theta_{2n}$;
therefore, this number is an \emph{invariant} of the line bundle
$\lambda$.

\begin{conjecture}
$\Tr_\lambda(\boldsymbol1_\lambda)=C(M)\cdot\chi(\lambda)$ where $C(M)$ does not
depend on $\lambda$ and $\chi(\lambda)$ is the Euler characteristic
of $\lambda$.
\end{conjecture}
This Conjecture is based on the discussions with Boris Feigin.

\section{Computation in the case $M=\CP^n$}

Unfortunately, the author has not found a simple computation of the
integrals $\int_{\langle\C P^n\rangle}\Theta_{2n}$ for projective spaces by the
method of the Dolbeault complex, described in Section 3. On the
other hand, the computation using the \v Cech complex turns
out to be relatively simple. Therefore, our viewpoint is a
compromise: we don't develop the general theory of the
integration by the method of the \v Cech complex, but in the case
of the projective spaces we reprove the analogs of Theorem 2.3
and Corollary 2.3 for the \v Cech method.

The higher cohomology of the sheaf of holomorphic differential
operators on $\CP^n$ vanishes, and we work with algebra of global
holomorphic differential operators.

\subsection{The explicit construction of the integral}

We will work with the covering
\begin{multline*}
U_1\cup U_2\cup\ldots\cup U_{n+1}=\CP^n=\{\,(z_1,\dots,z_{n+1})
/{\sim}\mid z_i\in\C,\\
(i=1,\dots,n+1),\quad\text{not all
$z_i$ are equal to~$0$}\,\}.
\end{multline*}
where
$$
U_i=\{\,(z_1,\dots,z_{n+1})\mid z_i\ne0\,\}
$$
For all $i=1,\dots,n+1$ and all the points $z^0\in U_i$ we choose
a holomorphic coordinate system $\iota_{z^0,U_i}\colon
U_i\to\C^n$ ($\iota_{z^0,U_i}(z^0)=0$), depending holomorphically
on the point~$z^0$ (for fixed~$U_i$).

We set
\begin{equation}
\iota_{z^0,U_i}(z_1,\dots,z_{n+1})=\left(\,\nfrac{z_1}{z_i}-
\nfrac{z_1^0}{z_i^0},\dots,\nfrac{z_{i-1}}{z_i}-\nfrac{z_{i-1}^0}{z_i^0},
\nfrac{z_{i+1}}{z_i}-\nfrac{z_{i+1}^0}{z_i^0},\dots,
\nfrac{z_{n+1}}{z_i}-\nfrac{z_{n+1}^0}{z_i^0}\,\right).
\end{equation}

We define the cocycles
$$
\iota^*_{z^0,U_{j_1}}(\xi),\dots,\iota^*_{z^0,U_{jk}}(\xi)\in
C_\Lie^l(\gl_\infty^\fin(\Dif_{\CP^n});\C)
$$
for any cocycle $\xi\in C_\Lie^l(\gl_\infty^\fin(\Dif_n);\C)$ and
any point $z^0\in U_{j_1}\cap\ldots\cap U_{j_k}$.

In our choice, $\iota_{z^1,U_j}=\iota_{z^0,U_i}\circ A$ where
$A$ is a global holomorphic automorphism of~$\CP^n$,
$A\in\GL_{n+1}(\C)/\C^*$, for any $i,j$ and $z^0\in U_i$, $z^1\in
U_j$.

\emph{Construction of the integral}:

I. For any $z^0\in U_1\cap U_i$ ($i>1$) we find a cochain
$$
\Theta_{1,i}(z^0)\in C_\Lie^{l-1}
(\gl^\fin_\infty(\Dif_{\CP^n});\C)
$$
such that
\begin{equation}
\iota^*_{z^0,U_i}\xi-\iota^*_{z^0,U_1}\xi=\delta_\Lie\Theta_{1,i}(z^0).
\end{equation}
Furthermore, for $z_0\in U_1\cap U_2\cap U_i$ ($i>2$)
we find a cochain
$$
\Theta_{1,2,i}(z^0)\in
C_\Lie^{l-2}(\gl_\infty^\fin(\Dif_{\CP^n});\C)
$$
such that
\begin{gather}
\Theta_{1,i}(z^0)-\Theta_{1,2}(z^0)=\delta_\Lie\Theta_{1,2,i}(z^0)\\
\hbox to 5cm{\dotfill}\nonumber
\end{gather}

Finally, for $z^0\in U_1\cap\ldots\cap U_{n+1}$ we find
$$
\Theta_{1,2,\dots,n+1}(z^0)\in
C_\Lie^{l-n}(\gl_\infty^\fin(\Dif_{\CP^n});\C)
$$
such that
\begin{equation}
\Theta_{1,2,\dots,n-1,n+1}(z^0)-\Theta_{1,2,\dots,n-1,n}(z^0)=
\delta\Theta_{1,2,\dots,n+1}(z^0).
\end{equation}

II. For any point $z_0\in U_1\cap U_2\cap\ldots\cap U_{n+1}$ we
find a solution
$\{\Theta^1_{1,2,\dots,n+1},\Theta^2_{1,2,\dots,n+1},\dots,
\Theta^n_{1,2,\dots,n+1}\}$,
$$
\Theta^i_{1,2,\dots,n+1}\in
C_\Lie^{l-n-i}(\gl_\infty^\fin(\Dif_{\CP^n});\C)\otimes_\C
\Omega^i_{U_1\cap\ldots\cap U_{n+1}}
$$
of the following system:
\begin{equation}
\begin{cases}
d_{\DR}\Theta_{1,2,\dots,n+1}=\delta_\Lie\Theta^1_{1,2,\dots,n+1}\\
d_{\DR}\Theta^1_{1,2,\dots,n+1}=\delta_\Lie\Theta^2_{1,2,\dots,n+1}\\
\hbox to 5cm{\dotfill}\\
d_\DR\Theta_{1,2,\dots,n+1}^{n-1}=\delta_\Lie\Theta^n_{1,2,\dots,n+1}
\end{cases}
\end{equation}
(here we denote by $\Omega^i_{U_1\cap\ldots\cap U_{n+1}}$ the space
of \emph{holomorphic} $i$-forms on $U_1\cap\ldots\cap U_{n+1}$).

We set $T^n=\{\,z_1=1,\ |z_2|=\ldots=|z_{n+1}|=1\,\}\subset\CP^n$.
We want to find conditions for which the following two statements hold:

(i) $\int_{T_n}\Theta^n_{1,2,\dots,n+1}\in C_\Lie^{l-2n}
(\gl_\infty^\fin(\Dif_{\CP^n});\C)$ is a \emph{cocycle};

(ii) if $[\xi]=0$ then $\left[\int_{T^n}\Theta^n_{1,2,\dots,n}\right]=0$
(here $[\dots]$ denotes the cohomological class of the cocycle).

We find these conditions and prove both statements in the next Subsection.

\subsection{Main Theorem}

The following reformulation of the construction
of the integral due on the previous Subsection will be very useful in
the course of the proof of statements (i) and~(ii).

Let $\sigma^i=\{\,(x_1,\dots,x_{i+1})\in\R^{i+1}\mid
\sum x_1+x_2+\ldots+x_{i+1}=1\,\}$ be the $i$-simplex, we denote by
$[1],[2],\dots,[i+1]$ its vertices, by $[1,2]$, $[1,3]$, $\dots$ ---
its $1$-faces, and so on.

Let us define the space $\wCP^n$ by the following way. The space~$\wCP^n$
is glued from the following pieces:
\begin{align*}\relax
&U_1,\dots,U_{n+1};\\
&(U_i\cap U_j)\times\sigma^1\qquad(i<j);\\
&(U_i\cap U_j\cap U_k)\times\sigma^2\qquad(i<j<k);\\
&\hbox to 5cm{\dotfill}\\
&(U_1\cap U_2\cap\ldots\cap U_{n+1})\times\sigma^n.
\end{align*}
We glue to the space
$$
(U_1\cap\ldots\cap U_{n+1})\times[1,2,\dots,n]
\subset(U_1\cap\ldots\cap U_{n+1})\times\sigma^n
$$
(here $[1,2,\dots,n]$ is a $(n-1)$-face of the simplex~$\sigma^n$)
the ``piece'' $(U_1\cap\ldots\cap U_n)\times\sigma^{n-1}$,
to the space $(U_1\cap\ldots\cap U_{n+1})\times[1,2,\dots,n-1,n+1]$
the ``piece'' $(U_1\cap\ldots\cap U_{n-1}\cap U_{n+1})\times\sigma^{n-1}$,
and so on. We glue all the pieces to each other along the faces
of simplexes, and after that we obtain a \emph{connected} space, we
denote it by~$\wCP^n$.

The projection $p\colon\wCP^n\to\CP^n$ is well-defined
(the fibers are simplexes of different dimension).

We connect a holomorphic coordinate system in the neighbourhood of
the point $p(\wz)$ with any point $\wz\in\CP^n$ in the following way.

First of all, let us suppose, that
$\wz\in(U_1\cap\ldots\cap U_{n+1})\times\sigma^n$, $z^0=p(\wz)$.
The coordinate systems, connected with points
$z^0\times[1]$, $z^0\times[2]$, $\dots$, $z^0\times[n+1]$
are the coordinate systems
$\iota_{z^0,U_1},\dots,\iota_{z^0,U_{n+1}}$ correspondently.

Furthermore, we define a coordinate systems on whole simplex
$z^0\times\sigma^n$ such that the following conditions (1)--(4) hold:

(1) the coordinate system depends smoothly on the point
of~$\sigma^n$;

(2) the coordinate systems, connected with the points
$z^0\times[1]$, $\dots$, $z^0\times[n+1]$, are
$\iota_{z^0,U_1},\dots,\iota_{z^0,U_{n+1}}$;

(3) for fixed $t\in\sigma^n$,
the coordinate system in the point $z^0\times t$ depends
\emph{holomorphically} on~$z^0$;

(4) if $t\in[\iota_1,\dots,\iota_k]\subset\sigma^n$,
then coordinate system on $(U_1\cap\ldots\cap U_{n+1})\times t$
can be extended to a \emph{holomorphic} coordinate system on
$(U_{i_1}\cap\ldots\cap U_{i_k})\times t$.

Let us suppose, that conditions (1)--(4) hold. Then one can associate to
any point~$\wz$ a holomorphic coordinate system in the neighbourhood
of the point~$p(\wz)$. Therefore, for a fixed cocycle
$\xi\in C_\Lie^l(\gl_\infty^\fin(\Dif_n);\C)$
one can associate the cocycle
$\xi(\wz)\in C_\Lie^l(\gl_\infty^\fin(\Dif_{\CP^n});\C)$.

We consider the \emph{canonical} solution $\{\Theta_1,\dots,\Theta_{2n}\}$,
$\Theta_i\in\Omega^i_{\wCP^n}\otimes C_\Lie^l(\gl_\infty^\fin
(\Dif_{\CP^n});\C)$, of the system:
\begin{equation}
\begin{cases}
d_\DR\xi(\wz)=\delta_\Lie\Theta_1\\
d_\DR\Theta_1=\delta_\Lie\Theta_2\\
\hbox to 3cm{\dotfill}\\
d_\DR\Theta_{2n-1}=\delta_\Lie\Theta_{2n}
\end{cases}
\end{equation}
(see formula (14)).

The $2n$-form $\Theta_{2n}$ on $\wCP^n$ depends \emph{holomorphically}
on~$p(\wz)$, therefore, $\Theta_{2n}=0$ for
$\wz\not\in(U_1\cap\ldots\cap U_{n+1})\times\sigma^n$.

\begin{theorem}
If conditions \emph{(1)--(4)} above hold, then\emph:

\emph{(i)} $\int_{T^n\times\sigma^n}\Theta_{2n}\in C^{l-2n}_\Lie
(\gl_\infty^\fin(\Dif_{\CP^n});\C)$
is a cocycle\emph;

\emph{(ii)} if $[\xi]=0$ then $\left[\int_{T^n\times\sigma^n}
\Theta_{2n}\right]=0$.
\end{theorem}

(Here $T^n\subset U_1\cap\ldots\cap U_{n+1}$,
$T^n=\{(1,z_2,\dots,z_{n+1})\mid|z_2|=\ldots=|z_{n+1}|=1$.)

\begin{proof}
(ii) is a direct consequence of~(i): one can repeat the arguments of
Theorem~2.3.(1). Let us prove~(i).

We have:
$$
\delta_\Lie\left(\int_{T^n\times\sigma^n}\Theta_{2n}\right)=
\int_{T^n\times\sigma^n}\delta_\Lie\Theta_{2n}=
\int_{T^n\times\sigma^n}d_\DR\Theta_{2n-1}=
\int_{T^n\times\partial(\sigma^n)}\Theta_{2n-1}.
$$
If a point~$t$ lies in a $(n-1)$-face of the simplex~$\sigma^n$,
say $t\in[1,2,\dots,n]$, then one can extended $\Theta_{2n-1}$ from
$(U_1\cap\ldots\cap U_{n+1})\times[1,2,\dots,n]$ to a $(2n-1)$-form
on \hbox{$(U_1\cap\ldots\cap U_n)\times[1,2,\dots,n]$}, which depends
\emph{holomorphically} on the first argument.
Therefore, the $(2n-1)$-form $\Theta_{2n-1}$ is defined on
$(S^1\times\ldots\times S^1\times D^2)\times[1,2,\dots,n]$
and depends \emph{holomorphically} on the point of $z\in D^2$
via condition~(4) (here $D^2=\{\,z\in\C\mid|z|\le1\,\}$,
$\partial D^2=S^1$). Therefore
$\int_{T^n\times[1,2,\dots,n]}\Theta_{2n-1}=0$
according to the Cauchy Theorem. The proof is the same for other
$(n-1)$-faces in $\partial\sigma^n$.
\end{proof}

In the notations of Subsec.~4.1 we have:
\begin{align*}\relax
&\int_{z^0\times[1,i]}\Theta_1=\Theta_{1,i}(z_0),\\
&\int_{z^0\times[1,2,i]}\Theta_2=\Theta_{1,2,i}(z^0),\\
&\hbox to 3cm{\dotfill}\\
&\int_{z^0\times[1,2,\dots,n+1]}\Theta_n=\Theta_{1,2,\dots,n}(z_0).
\end{align*}
Therefore, the definition of the integral, given in this Subsection,
coincides with the definition, given in Subsec.~4.1.

\begin{remark}
The Theorem remains true for the Lie algebra
$\gl_\infty^\fin(\Dif_{\lambda,\CP^n})$
for any holomorphic line bundle~$\lambda$ on~$\CP^n$.
The proof is the same.
\end{remark}

\begin{lemma}
The coordinate systems
$\iota_{z^0,U_1},\dots,\iota_{z^0,U_{n+1}}$, defined by the
formula~\emph{(22)},
satisfy conditions \emph{(1)--(4)} above.
\end{lemma}

\begin{proof}
Straightforward.
\end{proof}

\subsection{Computation for $\CP^1$}

\subsubsection{}
We denote $U_1$ by $V$ and $U_2$ by~$U$, and let
$z^0=(x_0,y_0)$ be a point $U\cap V\subset\CP^1$.
Let us find $\Theta_{1,2}(x_0,y_0)$ (see Subsection~4.1).

If $\iota_{z^0,U},\iota_{z^0,V}$ are the coordinate systems,
defined by~(22), then the maps
$j_{z^0,U},j_{z^0,V}\colon\Dif_{\CP^1}\to\Dif_1$
are defined by the formulas:
\begin{equation}
\begin{gathered}
j_{z^0,U}(\D)(f)=(\iota^{-1}_{z^0,U})^*\D(\iota^*_{z^0,U}f),\\
j_{z^0,V}(\D)(f)=(\iota^{-1}_{z^0,V})^*\D(\iota^*_{z^0,V}f).
\end{gathered}
\end{equation}

We have: $\iota_{z^0,V}=\iota_{z^0,U}\circ A_{z^0}$,
where $A\in\GL_2(\C)/\C^*$ is the global automorphism
of~$\CP^1$, defined by the matrix:
$$
A_{z^0}=
\begin{pmatrix}
\nfrac{x_0^2}{y_0^2}-1&\nfrac{x_0}{y_0}\\
\nfrac{x_0}{y_0}&0
\end{pmatrix}.
$$
From (28) we have:
\begin{equation}
j_{z^0,V}(\D)=j_{z^0,U}(A_{z^0}^{-1}\D A_{z^0}^{-1})
\end{equation}
Here the differential operator $A^{-1}_{z^0}\D A_{z^0}$
defined by a differential operator~$\D$ by the formulas:
$A_{z^0}^{-1}\D A_{z^0}(f)=\wf$, where
\begin{equation}
\wf(x)=\hf(A_{z^0}^{-1}x)\quad\text{and}\quad\hf(x)=\D(f\circ A_{z^0}(x)).
\end{equation}
Denote $A^{-1}\D A$ by $\Ad(A^{-1})(\D)$.

We have:
$$
\iota^*_{z^0,U}\xi-\iota^*_{z^0,V}\xi=
\iota^*_{z^0,U}\xi-\Ad(A^{-1}_{z^0})\cdot\iota^*_{z^0,U}\xi.
$$

\subsubsection{} Let $A_t$ be a matrix \emph{piece-wise}
smooth path, $A_t\in\GL_n(\C)$, $t\in[0,1]$. We have:
\begin{multline}
(\Ad(A_1)-\Ad(A_0))(x)=
\sum_{\dots}(A_{t+\epsilon}xA_{t+\epsilon}^{-1}-A_txA_t^{-1})=\\
=\sum_{\dots}(\Ad(A_{t+\epsilon}A_t^{-1})-1)(A_tx A_t^{-1})=
\int_{t\in[0,1]}\ad(a_t)(A_txA_t^{-1})
\end{multline}
where
\begin{equation}
a_t=\nfrac d{d\epsilon}(A_{t+\epsilon}A_t^{-1})\big|_{\epsilon=0}.
\end{equation}

Let us choose a path, connecting $\Id$ and~$A_{z^0}$.

Suppose:
\begin{align*}\relax
&X_t=
\begin{pmatrix}
e^{\pi it}&0\\
0&1
\end{pmatrix}\ ,\qquad t\in[0,1],\\
&Y_t=
\begin{pmatrix}
-\cos\nfrac\pi2t&\sin\nfrac\pi2t\mph\\
\sin\nfrac\pi2t&\cos\nfrac\pi2t\mph
\end{pmatrix}\ ,\qquad t\in[0,1],\\
&Z_t=
\begin{pmatrix}
t^2-1&t\\
t&0
\end{pmatrix}\ ,\qquad
\parbox{5cm}{$t$ be a complex path, connecting $1$ and~$\nfrac{x_0}{y_0}$.}
\end{align*}
We have:
$$
X_0=\Id,\quad X_1=Y_0,\quad Y_1=Z_1,\quad Z_{\frac{x_0}{y_0}}=A_{z^0}.
$$
According to the formulas (31) and (32), we have:
\begin{multline}
\iota^*_{z^0,U}\xi-\iota^*_{z^0,V}\xi=
-\int_0^1\ad(x_t)(\Ad(X_t)^{-1}(\iota^*_{z^0,U}\xi))-\\
-\int_0^1\ad(y_t)(\Ad(Y_t^{-1})(\iota^*_{z^0,U}\xi))-
\int_1^{x_0/y_0}\ad(z_t)(\Ad(Z_t^{-1})(\iota^*_{z^0,U}\xi)).
\end{multline}
In the affine coordinate $z=\nfrac xy$ on~$U$ we have:
\begin{align*}\relax
&x_t=-\pi i z\pp{}z,\\
&y_t=\nfrac\pi2z\pp{}z,\\
&z_t=-\left(1+\nfrac1{t^2}\right)z^2\pp{}z.
\end{align*}

Let us define the $(l-1)$-cochain
\begin{equation}
\phi_X(\D_1,\dots,\D_{l-1})=
\xi\left(j_{z^0,U}(X_t^{-1}x_tX_t),j_{z^0,U}(X_t^{-1}\D_1X_t),
\dots,j_U(X_t^{-1}\D_{l-1}X_t)\right).
\end{equation}
Then
\begin{equation}
\ad(x_t)(\Ad(X_t^{-1})(\iota^*_{z^0,U}))(\D_1,\dots,\D_{l-1})=
(\delta_\Lie\phi_X)(\D_1,\dots,\D_{l-1}).
\end{equation}
There are analogous formulas for $Y_t$ and~$Z_t$.

As we will see below, the values of $\int_{S^1}\Theta^1_{1,2}$
depend only on the $Z_t$-term.

\subsubsection{}
Denote $S^1=\left\{\,(x,y)\in\CP^1\,\Big|\,\left|\nfrac xy\right|=1\,\right\}$.
We consider the coordinate systems $\iota_{z^0,U}$ in all the points
$z^0\in S^1$. The (global) holomorphic vector field on~$\CP^1$,
corresponding to a tangent vector to~$S^1$ in any point (in the
sense of Subsec.~2.1) is equal, up to a constant factor, $\pp{}z$
(in the affine coordinate in~$U$).

Therefore, the $X_t$- and $Y_t$-terms in~$\Theta_{1,2}^1$
have, up to a constant factor, the following form:
$$
\int\Ad(X_t^{-1})\xi\left(j_{z^0,U}\left(z\pp{}z\right),
j_{z^0,U}\left(\pp{}z\right),j_{z^0,U}(\D_1),\dots,j_{z^0,U}(\D_{l-2})
\right)
$$
and
$$
\int\Ad(Y_t^{-1})\xi\left(j_{z^0,U}\left(z\pp{}z\right),
j_{z^0,U}\left(\pp{}z\right),j_{z^0,U}(\D_1),\dots,j_{z^0,U}(\D_{l-2})
\right).
$$
These terms have nonzero grading as expression of the form
$\Psi(\D_1,\dots,\D_{l-2})$,
and therefore the sequel computations do not depend on these terms.
$Z_t$-term has a form
$$
\int\Ad(Z_t^{-1})\xi\left(j_{z^0,U}\left(z^2\pp{}z\right),
j_{z^0,U}\left(\pp{}z\right),j_{z^0,U}(\D_1),\dots,
j_{z^0,U}\left(\D_{l-2}\right)\right)
$$
and this expression has the grading 0.

Let us calculate $Z_t$-term. We have:
$$
(Z_t)=\int_{\nfrac {x_0}{y_0}\in
S^1}\int^{\varphi_0}\Ad(Z^{-1}_t)\xi\left(j_{z^0,U}\left(\pp{}z\right),
j_{z^0,U}(Z_t),j_{z^0,U}(\D_1),
\dots,j_{z^0,U}(\D_{l-2})\right)
$$
here $\nfrac{x_0}{y_0}=e^{2\pi\varphi_0}$, $\varphi_0\in [0,1]$,
$Z_t=-\left(1+\nfrac{1}{t^2}\right)z^2\pp{}z$.

\emph{Calculation of the integral}:
\begin{multline*}
(Z_t)=\int_{\varphi_0\in[0,1]}\int^{\varphi_0}_0\left(1+\nfrac{1}{e^{4\pi
i\beta}}\right)\cdot\\
\left(\pp{}z,z^2\pp{}z,j_{z^0,U}(\Ad(Z^{-1}_t)(\D_1)),\dots,j_{z^0,U
}
(\Ad(Z^{-1}_t)(\D_{l-2}))\right)\,d\beta d\varphi_0.
\end{multline*}
(here $t=e^{2\pi i\beta}$).

\subsubsection{} Let $\xi=\Psi_{2k+1}$ $(k\ge 1)$ be the Lifting
cocycle on the Lie algebra $\gl^\fin_\infty(\Dif_1)$ (see
Subsec.~1.2), and
$\D_1,\dots,\D_{2k-1}\in\gl^\fin_\infty\otimes 1$. Then
$\Ad(Z^{-1}_t)(\D_i)=\D_i$, and we are able to conclude the
calculation. We have:
\begin{equation}
(Z_t)=\left(\nfrac{1}{2}+\nfrac{1}{4\pi
i}\right)\Psi_{2k+1}\left(\pp{}z, z^2\pp{}z,\D_1,\dots,\D_{2k-1}\right)
\end{equation}
where $\D_i\in\gl^\fin_\infty\otimes 1$.

\emph{It follows directly from our definitions that $\pp{}z$ and
$z^2\pp{}z$ in~\emph{(36)} are infinite matrices $\Id\otimes\pp{}z$ and
$\Id\otimes z^2\pp{}z$. Strictly speaking, these matrices do not lie
in the Lie algebra $\gl^\fin_\infty(\Dif_1)$ \emph(see also}~(21)).

\begin{theorem} The cochain
$\Psi_{2k+1}\left(\pp{}z, z^2\pp{}z,\D_1\dots,\D_{2k-1}\right)$ has nonzero value
on some \hbox{$(2k-1)$-cycle} on the Lie algebra
$\gl^\fin_\infty\otimes 1$.
\end{theorem}

\begin{proof} First of all we consider the simplest case $k=1$.
We have to calculate $\Psi_3(\partial, x^2\partial,1)$ (here
$\partial=\partial\otimes\Id$,
$x^2\partial=x^2\partial\otimes\Id$, $1=E_{11}$. We have:
$$
\Psi_3(A_1,A_2,A_3)=\Altl_{A,D} \Tr(D_1A_1\cdot D_2A_2\cdot
A_3)+\Altl_A \Tr(Q\cdot A_1\cdot A_2\cdot A_3),
$$
where $D_1=\ad\ln\partial$, $D_2=\ad\ln x$,
$$
Q=x^{-1}\partial^{-1}+\nfrac{1}{2}x^{-2}\partial^{-2}+\ldots+
\nfrac{(n-1)!}{n}x^{-n}\partial^{-n}+\dots
\qquad\text{(see Subsec.~1.2).}
$$

The first summand in $\Psi_3(\partial,x^2\partial,1)$ is equal
to~$-2$, the second one is equal to~$-1$. Therefore,
$\Psi_3(\partial,x^2\partial,1)=-3\ne0$.

The general case:

let $\suml_iA_1^{(i)}\wedge\dots\wedge A^{(i)}_{2k-1}$ be a cycle
in $\gl_\infty^\fin\otimes1$, we have to calculate
$$
\sum_i\Psi_{2k+1}(\partial,x^2\partial,A_1^{(i)},\dots,A^{(i)}_{2k-1})
\qquad\text{(see Subsec.~1.2).}
$$
The ``leading'' term $\Altl_{A,D}\Tr(D_1A_1\cdot D_2A_2\cdot
A_3\cdot\ldots\cdot A_{2k+1})$ is equal to
$-2\Tr\left(\suml_i\Altl(A_1^{(i)}\cdot A_2^{(i)}\cdot\ldots\cdot
A_{2k-1}^{(i)})\right)$. Furthermore, we consider the term
$$
\Altl_{A,D}(D_1A_1\cdot A_2\cdot\ldots\cdot A_{2s-1}\cdot
D_2A_{2s}\cdot A_{2s+1}\cdot\ldots\cdot A_{2k+1}).
$$
It is necessarily that $A_1=x^2\partial\otimes\Id$,\ \
$A_{2s}=\partial\otimes\Id$ (otherwise this expression is equal
to~$0$), and matrices $D_1(x^2\partial\otimes\Id)$ and
$D_2(\partial\otimes\Id)$ commute with $\gl_\infty^\fin\otimes1$.
Therefore, every such term is equal to
$$
-2\Tr\left(\sum_i\Altl(A_1^{(i)}\cdot\ldots\cdot
A_{2k-1}^{(i)})\right).
$$
The remaining term $\Altl_A\Tr(Q\cdot A_1\cdot
A_2\cdot\ldots\cdot A_{2k+1})$ is equal to~$-1$. Therefore,
{\multlinegap=0pt
\begin{multline*}
\sum_i\Psi_{2k+1}(\partial\otimes\Id,x^2\partial\otimes\Id,A_1^{(i)},
\dots,A_{2k-1}^{(i)})=\\
=\nfrac12(k+1)\cdot\left(-2\Tr\left(\sum_i\Altl(A_1^{(i)}\cdot\ldots
\cdot A_{2k-1}^{(i)})\right)\right)-\Tr\left(\sum_i\Altl(A_1^{(i)}
\cdot\ldots\cdot A_{2k-1}^{(i)})\right)=\\
=-(k+2)\cdot\Tr\left(\sum_i\Altl(A_1^{(i)}\cdot\ldots\cdot
A_{2k-1}^{(i)})\right).
\end{multline*}}%
It follows from the standard facts on the cohomology of the Lie
algebra~$\gl_\infty^\fin$ (see~[F], Ch.~2.1) that the last
expression is not equal to~$0$ for some cycle
$\suml_iA_1^{(i)}\wedge\dots\wedge A_{2k-1}^{(i)}$.
\end{proof}

\subsubsection{}

Using the same methods, one can show, that the value of
$$
\int_{\CP^1}\Psi_{2k+1}\in
C_\Lie^{2k-1}(\gl_\infty^\fin(\Dif_{\lambda,\CP^1});\C)
$$ 
on the
matrix cycle $\suml_iA_1^{(i)}\wedge\dots\wedge A_{2k-1}^{(i)}$
(here $\lambda=\O(\nl)$) is equal, up to a nonzero undepending
on~$\nl$ constant, to
\begin{multline*}
\sum_i\Psi_{2k+1}(\Id\otimes\partial,\Id\otimes(x^2\partial-\nl
x),A_1^{(i)},\dots,A_{2k-1}^{(i)})=\\
=-(k+2)(\nl+1)\cdot\left(\sum_i\Tr\Altl A_1^{(i)}\cdot\ldots\cdot
A_{2k-1}^{(i)}\right).
\end{multline*}
Let us note, that the number $\nl+1$ is equal to the Euler
characteristic $\chi(\O(\nl))$.

\subsubsection{}

\begin{theorem}
$$
H^\ndot(\gl_\infty^\fin(\Dif_1);\C)=\wedge^\ndot(\Psi_3,\Psi_5,\Psi_7,\dots).
$$
\end{theorem}

\begin{proof}
It follows from Theorem 4.2 and Theorem~4.3.4 that the cocycles
$\Psi_3,\Psi_5,\Psi_7,\dots$ are noncohomologous to zero cocycles
on the Lie algebra $\gl_\infty^\fin(\Dif_1)$.

Furthermore, the standard \emph{homomorphism of the Lie algebras}
$$
\gl_\infty^\fin(\Dif_1)\oplus\gl_\infty^\fin(\Dif_1)\to
\gl_\infty^\fin(\Dif_1)
$$ 
defines the Hopf algebra structure on the
cohomology $H^\ndot(\gl_\infty^\fin(\Dif_1);\C)$. We denote
by~$\Delta$ the coproduct in this Hopf algebra.

\begin{lemma}
The Lifting cocycles $\Psi_3,\Psi_5,\Psi_7,\dots$ are primitive
elements in $H^\ndot(\gl_\infty^\fin(\Dif_1);\C)$, i.\,e.\
$\Delta\Psi_{2k+1}=1\otimes\Psi_{2k+1}+\Psi_{2k+1}\otimes1$.
\end{lemma}

\begin{proof}
It is obvious.
\end{proof}

It is well-known theorem, that any (super-) commutative and
cocommutative Hopf algebra~$A$ is equal to (super-) symmetric
algebra, generating by its primitive elements (i.\,e.\ such
elements $a\in A$ that $\Delta a=1\otimes a+a\otimes1$).
Therefore, there do not exist any nontrivial relations between
$\Psi_3,\Psi_5,\Psi_7,\dots$ in
$H^\ndot(\gl_\infty^\fin(\Dif_1);\C)$. On the other hand,
according to~[FT1] $H^\ndot(\gl_\infty^\fin(\Dif_1);\C)$ is the
exterior algebra with the unique generator in each dimension
$3,5,7\dots$. It proves, that
$$
H^\ndot(\gl_\infty^\fin(\Dif_1);\C)=\wedge^\ndot(\Psi_3.\Psi_5,\Psi_7\dots).
$$
\end{proof}

\subsubsection{}

\begin{theorem}
Let $\iota\colon\Dif_{\CP^1}\to\Dif_1$ is the map, connected with
a choice of coordinate system in any point of~$\CP^1$. Then
$$
H^\ndot(\gl_\infty^\fin(\Dif_{\CP^1});\C)\simeq
\wedge^\ndot\left(\int_{\CP^1}\Psi_3;\ \iota^*\Psi_3,
\int_{\CP^1}\Psi_5;\ \iota^*\Psi_5,\int_{\CP^1}\Psi_7;\ \dots\right)
$$
\end{theorem}

\begin{proof}
First of all, one can show by the method of~[FT1], that
$H^\ndot(\gl^\fin_\infty(\Dif_{\CP^1});\C)$ is the exterior
algebra with the unique generator is dimension~$1$ and two
generators in each dimension $3,5,7,\dots$. Furthermore, there
exists the standard map of the Lie algebras
$\gl_\infty^\fin(\Dif_{\CP^1})\oplus\gl_\infty^\fin(\Dif_{\CP^1})\to
\gl_\infty^\fin(\Dif_{\CP^1})$, which defines the Hopf algebra
structure on cohomology
$H^\ndot(\gl^\fin_\infty(\Dif_{\CP^1});\C)$. It is easy to prove
that $\int_{\CP^1}\Psi_{2k+1}$ and $\iota^*\Psi_{2k+1}$ ($k\ge1$)
are \emph{primitive} elements in
$H^\ndot(\gl^\fin_\infty(\Dif_{\CP^1});\C)$ with respect to the Hopf
algebra structure. It follows from Theorem~4.3.4, that
$\int_{\CP^1}\Psi_{2k+1}$ ($k\ge1$) are nonzero elements in
$H^\ndot(\gl_\infty^\fin(\Dif_{\CP^1});\C)$. It remains to prove that
$[\iota^*\Psi_{2k-1}]\ne0$, $k\ge2$, and that the cocycles
$\iota^*\Psi_{2k-1}$ and $\int_{\CP^1}\Psi_{2k+1}$ are not
cohomologous.

To prove the first statement, let us note, that the integral
$\int_{\CP^1}\Psi_{2k-1}$ is defined via the pull-backs
$\iota^*\Psi_{2k-1}$, connected with all the points of~$\CP^1$.
It follows from Lemma~2.1 that all these cocycles (connected to
different points and a different choice of coordinate systems) are
cohomologous (Strictly speaking, we have to use the special class
of admissible coordinate systems, $\dots$ but is sufficient for
the definition of the integral.). Therefore, it follows from
Theorem~2.3 and Corollary~2.3 that if the cocycle
$\iota^*\Psi_{2k-1}$ is cohomologous to~$0$ than
$\int_{\CP^1}\Psi_{2k-1}$ is also cohomologous to zero, with the
contradiction to Theorem~4.3.4.

To prove the second statement, let us note, that the value of
$\int_{\CP^1}\Psi_{2k+1}$ on some \emph{matrix} cycle is not
equal to~$0$, but is obvious, that the value
of~$\iota^*\Psi_{2k-1}$ on any matrix cocycle is equal to~$0$.
\end{proof}

It is true also, that for $\lambda=\O(\nl)$, $\nl\ne-1$ and
(in the sense of the algebras of \emph{twisted} differential
operators $\Dif_{\lambda,\CP^1}$, $\nl\in\C$),
$6\left(\nfrac{\nl}2\right)^2+6\left(\nfrac{\nl}2\right)+1\ne0$
than
$$
H^\ndot(\gl_\infty^\fin(\Dif_{\lambda,\CP^1});\C)\simeq\wedge^\ndot
\left(\int_{\CP^1}\Psi_3;\ \iota^*\Psi_3,\int_{\CP^1}\Psi_5;\ 
\iota^*\Psi_5,\int_{\CP^1}\Psi_7;\ \dots\right).
$$
It is an interesting exercise to understand what appears when
$\nl=-1$ or
$6\left(\nfrac{\nl}2\right)^2+6\left(\nfrac{\nl}2\right)+1=0$.

\subsection{Computation for $\CP^n$}

Let us denote by $\Dif_{\lambda,\CP^n}$ the associative algebra
of global holomorphic differential operators in the bundle
$\O(\nl)$ on~$\CP^n$.

Using the methods, analogous to methods of Subsec.~4.3, one can
prove, that the value of the integral
$$
\int_{\CP^n}\Psi_{2k+1}\in
C_\Lie^{2k-2n+1}(\gl_\infty^\fin(\Dif_{\lambda,\CP^n});\C)\qquad(k\ge
n)
$$
on a \emph{matrix} cycle
$
\suml_iA_1^{(i)}\wedge\dots\wedge A^{(i)}_{2k-2n+1}$
($A_j^{(i)}\in\gl_\infty^\fin\otimes1$) is equal to
\begin{multline}
\sum_i\Psi_{2k+1}\left(\partial_1,\dots,\partial_n,x_1
\left(\sum_{j=1}^nx_j\partial_j\right)-\nl
x_1,\dots,x_n\left(\sum_jx_j\partial_j\right)-\nl
x_n,\right.\\
\left.A_1^{(i)},\dots,A_{2k-2n+1}^{(i)}\right).
\end{multline}
In (37) $\partial_i$ and
$x_i\left(\suml_{j=1}^nx_j\partial_j\right)-\nl x_i$ denotes
$\Id\otimes\partial_i$ and
$\Id\otimes\left(x_i\left(\suml_{j=1}^nx_j\partial_j\right)-\nl
x_i\right)$ in particular, these are \emph{infinite} matrices.
According to Remark~4.2, it is sufficient to prove that
$\int_{\CP^n}\Psi_{2k+1}\in
C_\Lie^{2k-2n+1}(\gl_\infty^\fin(\Dif_{\lambda,\CP^n});\C)$ is
not cohomologous to~$0$ in order to prove that $\Psi_{2k+1}\in
C_\Lie^{2k+1}(\gl_\infty^\fin(\Dif_n);\C)$ is not cohomologous
to~$0$ ($k\ge n$).

The difference between the general case of~$\CP^n$, $n\ge1$ and
the case of~$\CP^1$ is that it is difficult to compute the whole
polynomial on~$\nl$, arising in the formula~(37),  or
even its value for $\nl=0$ (Conjecture 3.3 gives us the
expected answer.) It is sufficient for  our aims to prove that
this polynomial is not equal to~$0$; therefore, it is sufficient
to prove that its \emph{leading coefficient} is not equal to~$0$.
For the calculation of the leading coefficient on~$\nl$
in~(37), it is sufficient to calculate
\begin{multline}
\sum_i\Psi_{2k+1}(\partial_1,\dots,\partial_n,\nl
x_1,\dots,\nl
x_n,A_1^{(i)},\dots,A_{2k-2n+1})=\\
\nl^n\sum_i\Psi_{2k+1}
(\partial_1,\dots,\partial_n,x_1,\dots,x_n,A_1^{(i)},\dots,
A_{2k-2n+1}^{(i)}).
\end{multline}
It is sufficient to find a matrix $(2k-2n+1)$-cycle
$\suml_iA_1^{(i)}\wedge\dots\wedge A_{2k-2n+1}^{(i)}$ such that
$$
\sum_i\Psi_{2k+1}(\partial_1,\dots,\partial_n,x_1,\dots,x_n,
A_1^{(i)},\dots,A_{2k-2n+1}^{(i)})\ne0
$$
(here $\partial_i$ and $x_i$ denotes $\Id\otimes\partial_i$ and
$\Id\otimes x_i$ respectively) ($k\ge n$).

First of all, let us consider the simplest case $k=n$. The
leading term in $\Psi_{2n+1}(A_1,\dots,A_{2n+1})$ is equal to
\begin{equation}
\Altl_{A,D}(D_1A_1\cdot\ldots\cdot D_{2n}A_{2n}\cdot
A_{2n+1})\qquad\text{(see Subsec.~1.3).}
\end{equation}
We set:
\begin{gather*}
D_1=\ad\ln\partial_1,\quad D_2=\ad\ln x_1,\quad
D_3=\ad\ln\partial_2,\quad D_4=\ad\ln x_2,\ \dots,\\
A_1=x_1,\quad A_2=\partial_1,\quad A_3=x_2,\quad A_4=\partial_2,\
\dots,\quad A_{2n+1}=1.
\end{gather*}
Then the leading term (39) is equal to $(-1)^n\cdot(2n)!$.
Furthermore, let us consider arbitrary term in
$\Psi_{2n+1}(A_1,\dots,A_{2n+1})$ (connected with the interval of
length $2n-2$ with marked points --- see Subsec.~1.3). We may
remove all the blocks of the form $A_i\cdot Q_{j,j+1}\cdot
A_{i+1}$ in the left-hand side of the expression, its value will
not change. Thus, we have the expression:
\begin{multline}
\Altl_{A,D}\Tr(A_1\cdot Q_{12}\cdot A_2\cdot A_3\cdot Q_{34}\cdot
A_4\cdot A_5\cdot Q_{56}\cdot A_6\cdot\ldots\cdot\\
\cdot A_{2l-1}\cdot
Q_{2l-1,2l}\cdot A_{2l}\cdot D_{2l+1}A_{2l+1}\cdot\ldots\cdot
D_{2n}A_{2n}\cdot A_{2n+1})
\end{multline}
(a sign does not appear from such a permutation). We have:
\begin{multline}
(40)=\Altl_{A,D}\Tr(A_1\cdot A_2\cdot Q_{12}\cdot A_3\cdot
A_4\cdot Q_{34}\cdot\ldots\cdot A_{2l-1}\cdot A_{2l}\cdot\\
\cdot Q_{2l-1,2l}\cdot A_{2n+1}\cdot D_{2l+1}A_{2l+1}\cdot\ldots\cdot
D_{2n}A_{2n}).
\end{multline}
We may suppose, that $A_{2n+1}=1$ in~(41) (all other terms are
annihilate with each other), and if
$$
A_{2s-1}=
\begin{cases}
x_j\\
\partial_j
\end{cases}
$$
than
$$
A_{2s}=
\begin{cases}
\partial_j\\
x_j\\
\end{cases}
$$
It is easy to see, that
$(41)=(-1)^n\cdot(l!)^2\cdot(2n-2l)!\cdot2^l$. In particular, all
the summands in
$\Psi_{2n+1}(\partial_1,\dots,\partial_n,x_1,\dots,x_n,1)$
contributes a nonzero number, and all these numbers have the same
sign. It proves that the leading coefficient on~$\nl$ in~(38)
is not equal to~$0$.

The general case of the cocycle~$\Psi_{2m+1}$ ($m\ge n$) is
deduced to the previous computation.

It follows from these computations and Theorem~4.2 that cocycles
$\Psi_{2n+1}$, $\Psi_{2n+3},\Psi_{2n+5},\dots$ are noncohomologous
to $0$ cocycles in $C_\Lie^\ndot(\gl_\infty^\fin(\Dif_n);\C)$. It
was proved in~[FT1] that $H^\ndot(\gl_\infty^\fin(\Dif_n);\C)$ is
the exterior algebra with the unique generator in each dimension
$2n+1,2n+3,2n+5,\dots$. Using the above arguments concerning the
Hopf algebra structure on $H^\ndot(\gl_\infty^\fin(\Dif_n);\C)$
we prove the following theorem:

\begin{theorem}
$$
H^\ndot(\gl_\infty^\fin(\Dif_n);\C)\simeq
\wedge^\ndot(\Psi_{2n+1},\Psi_{2n+3},\Psi_{2n+5},\dots).
$$
\end{theorem}

\qed


\begin{thebibliography}{XXX}



\bibitem[A]{A} Adler M. \emph{On a trace functional for formal
pseudo-differential operators and the symplectic structure of
the Korteweg -- de Vries type equations}, Invent.\  Math.\ 50 (1979), 219--248.

\bibitem[C]{C} Connes A., \emph{Cohomologie cyclique et
functeurs~$\Ext^n$}, C.R.\ Acad. Sci.\ Paris, 296, Ser.~I, pp.~953--958.

\bibitem[F]{F} Feigin B.L., \emph{Lie algebra $\gl(\lambda)$ and
cohomology of a Lie algebra of differential operators}, Russian
Mathematical Surveys 43(2), 1988, 169--170.


\bibitem[FT1]{FT1} Feigin B.L., Tsygan B.L. \emph{Cohomology of
the Lie algebra of the generalized Jacobian matrices} (in
Russian), Funct.\ Anal.\ Appl.\ 17(2) (1983), 86--87.

\bibitem[FT2]{FT2} Feigin B.L., Tsygan B.L., \emph{Riemann-Roch
theorem and Lie algebra cohomology}~I, in: \emph{Proceedings of
the winter school on Geometry and Physics}, Srni, 9--16 January
1988.

\bibitem[Fu]{Fu} Fuks D.B., \emph{Cohomology of
Infinite-Dimensional Lie algebras}, Consultants Bureau, New York
and London, 1986.

\bibitem[GF]{GF} Gelfand I.M., Fuks D.B., \emph{Cohomology of the
Lie algebra of tangent vector fields on a smooth manifold} II
(Russian), Funk.\ Anal.\ Appl., 4(2), 1970, 23--32. English transl.:
in I.M.\,Gelfand Collected Papers, Vol.~3.

\bibitem[GM]{GM} Gelfand I.M., Mathieu O., \emph{On the cohomology
of the Lie algebra of Hamiltonian vector fields}, J.~Funct.\
Anal.\ 108 (1992), pp.~347--360.

\bibitem[K]{K} Kassel C., \emph{Le r\'esidu non commutatif
\emph[d' apr\`es M.~Wodzicki}], S\'eminaire Bourbaki, 1988--89,
N708, Ast\'erisque 177--178, 1989, p.~199--229.


\bibitem[KK]{KK} Kravchenko O.S., Khesin B.A. \emph{Central
extensions of the algebra of pseudo-differential symbols}, Funct.\
Anal.\ Appl. 25(2), (1991), 83--85.

\bibitem[KLR]{KLR} Khesin B., Lyubashenko V., Roger C.,
\emph{Extensions and contractions of the Lie algebra of
\hbox{$q$-pseudodifferential} symbols}, J.~of Func.\ Anal., vol.~143
(1997), pp.~55--97.

\bibitem[LQ]{LQ} Loday J.L., Quillen D.G., \emph{Homologie
cyclique et homologie d'algebres de Lie des matrices}, C.R.\
Acad.\ Sci.\ Paris, 296, Ser.~I, 1983, pp.~295--297.



\bibitem[NT]{NT} Nest R., Tsygan B., \emph{Algebraic Index
Theorem for Families}, Preprint Kobenhavns Universitet, 
Preprint Series, 1993, No.\ 28.


\bibitem[Sh1]{Sh1} Shoikhet B., \emph{Cohomology of the Lie
algebras of differential operators\emph: Lifting formulas}, Amer.
Math.\ Soc.\ Transl., (2), Vol.~185, 1998,
and preprint q--alg/9712007.

\bibitem[Sh2]{Sh2} Shoikhet B., \emph{Lifting Formulas} II, to
appear in MRL, preprint math.\ QA/9801116.


\bibitem[Sh3]{Sh3} Shoikhet B., \emph{Certain topics on the Lie
algebra $\gl(\lambda)$ representation theory}, to appear in
Journal of Math.\ Science, and preprint q--alg/9703029.

\bibitem[T]{T} Tsygan B.L., \emph{Homology of the matrix Lie algebras
 over
rings and Hochschild homology} (Russian), Russian Mathematical
Surveys, 38(2), 1983, pp.~217--218.



\bibitem[W]{W} Wodzicki M., \emph{Noncommutative residue.
Fundamentals}, in: \emph{$K$-theory, Arithmetic and Geometry},
Yu.~I.~Manin (ed.), Lect.\ Notes in Math.\ 1289.








\end{thebibliography}
\end{document}